\documentclass[11pt,a4paper]{article}
\usepackage{authblk}
\usepackage{hyperref}
\usepackage{graphicx} 
\usepackage{mathtools} 
\usepackage{amsfonts} 
\usepackage[usenames,dvipsnames]{xcolor}
\usepackage{amsthm}
\usepackage{amssymb}
\usepackage[hmargin={26mm,26mm},vmargin={30mm,35mm}]{geometry}
\usepackage{thmtools}
\usepackage{cancel}
\usepackage[ruled,vlined]{algorithm2e}
\usepackage{siunitx} 
\usepackage{subcaption}
\usepackage{mathrsfs} 

\newcommand{\coloneq}{\mathrel{{:=}}}
\newcommand{\norm}[2][]{\|#2\|_{#1}}

\newcommand{\calH}{\mathcal{H}}
\newcommand{\Poly}[2][]{\mathbb{P}_{#1}^{#2}}

\newcommand{\Real}{\mathbb{R}}

\newcommand{\GRAD}{\vec{\nabla}}
\newcommand{\DIV}{\vec{\nabla} \cdot}

\newcommand{\GRADs}{\tens{\nabla}_{\!\mathrm s}}
\newcommand{\GRADss}{\tens{\nabla}_{\!{\rm ss}}}

\newcommand{\st}{\; : \;}

\newcommand{\tr}{\operatorname{tr}}
\newcommand{\argmin}{\operatorname*{arg\,min}}

\newcommand{\Th}[1][h]{\mathcal{T}_{#1}}
\newcommand{\Fh}[1][h]{\mathcal{F}_{#1}}
\newcommand{\FT}[1][h]{\mathcal{F}_T}

\newcommand{\Mh}[1][h]{\mathcal{M}_{#1}}

\newcommand{\normal}{\vec{n}}

\newcommand{\VT}{\underline{V}_T^0}
\newcommand{\Vh}{\underline{V}_h^0}

\newcommand{\uT}{\underline{\vec{u}}_T}
\newcommand{\vT}{\underline{\vec{v}}_T}
\newcommand{\uphi}{\underline{\phi}_T}
\newcommand{\uchi}{\underline{\chi}_T}

\newcommand{\pT}[1]{\mathrm{p}_T^{#1}}

\renewcommand{\vec}[1]{\boldsymbol{#1}}
\newcommand{\tens}[1]{\boldsymbol{#1}}

\newcommand{\uvec}[1]{\underline{\vec{#1}}}

\newcommand{\vVT}{\underline{\vec V}_T^1}
\newcommand{\vVh}{\underline{\vec V}_h^1}
\newcommand{\vVhD}{\uvec{V}_{h,\mathrm{D}}^1}
\newcommand{\vVhDn}[1]{\uvec{V}_{h,\mathrm{D}}^{1,#1}}
\newcommand{\vVhZ}{\uvec{V}_{h,0}^1}

\newcommand{\vpT}[1]{\vec{\mathrm{p}}_T^{#1}}

\newcommand{\tET}[1]{\tens{\mathrm{E}}_{T}^{#1}}
\newcommand{\DT}[1]{\mathrm{D}_T^{#1}}

\theoremstyle{remark}
\newtheorem*{remark}{Remark}

\newcommand{\email}[1]{\href{mailto:#1}{#1}}

\newcommand{\tF}{t_{\mathrm{F}}}

\title{A Hybrid-High Order method for fracture modelling}
\author[1,2]{Alessandra Crippa}
\author[2]{Julien Coatl\'{e}ven}
\author[1]{Daniele A. Di Pietro}
\author[2]{Nicolas Guy}
\author[2]{Soleiman Yousef}

\affil[1]{%
  Université de Montpellier, IMAG, Montpellier, France,
  \email{alessandra.crippa@etu.umontpellier.fr}, %
  \email{daniele.di-pietro@umontpellier.fr}
}
\affil[2]{%
  IFP Energies Nouvelles, Rueil-Malmaison, France,
  \email{julien.coatleven@ifpen.fr}, %
  \email{nicolas.guy@ifpen.fr}, %
  \email{soleiman.yousef@ifpen.fr}
}

\begin{document}

\maketitle

\begin{abstract}
  In this work, we introduce a new Hybrid High-Order method for the numerical simulation of fracture propagation based on phase-field models.
  The proposed method supports general meshes made of polygonal/polyhedral elements, which provides great flexibility in mesh design and adaptation, and can accommodate large variations of both the displacement and damage variables thanks to the use of fully discontinuous spaces.
  The resolution of the corresponding algebraic problem is based on a staggered time stepping scheme which takes advantage of static condensation for each subproblem.
  We provide extensive numerical validation of the method on classical two-dimensional fracture propagation problems, including a comparison with a more standard finite element scheme.
  \smallskip\\
  \textbf{MSC 2020:} %
  65N30, 
  65N08, 
  74R10, 
  \smallskip\\
  \textbf{Key words:} Fracture modelling, phase-field model, Hybrid High-Order methods, polyhedral methods

\end{abstract}


\section{Introduction}

In this work, we introduce a new Hybrid High-Order (HHO) method for the numerical simulation of fracture propagation based on phase-field models.
The proposed method supports general meshes made of polygonal/polyhedral elements, which provides great flexibility in mesh design and adaptation, and can accommodate large variations of both the displacement and damage variables thanks to the use of fully discontinuous spaces.
The corresponding algebraic problem is solved using a staggered time stepping scheme with static condensation for each subproblem.
\smallskip

Cracking is one of the most common failure mechanisms in engineering materials such as concrete, ceramics, and composites, affecting structures like wind turbines, dams, ships, containers, and buildings. Since the 1960s, several models and numerical methods have been developed to simulate fracture mechanics, which can be broadly categorized into discrete and continuous.

\emph{Discrete models}~\cite{Ngo.Scordelis:67} are quite intuitive in that they represent cracks as explicit discontinuities in the material, allowing the displacement field to be discontinuous across the crack surface.
However, such discontinuities can be challenging to handle numerically as they introduce a strong sensitivity to the mesh.
This issue can be mitigated resorting to, e.g.,
enriched finite elements~\cite{Belytschko.Black:99,Moes.Dolbow.ea:99,Moes.Belytschko:02,Wells.Sluys:01,Wu.Li.ea:15},
remeshing~\cite{Ngo.Scordelis:67,Bouchard.Bay.ea:00,Bouchard.Bay.ea:03},
or meshless approaches~\cite{Belytschko.Lu.ea:94}.
Other difficulties in discrete models include modeling crack initiation, propagation, and branching, especially in three-dimensional settings or complex geometries.

\emph{Continuous (smeared) crack models}~\cite{Rashid:68}, on the other hand, represent cracks as diffused zones within the material. The displacement jumps are smeared over a region of small but finite width, while material behaviour gradually deteriorates near the crack to model the degradation process.
This approach avoids the complexities associated with tracking discrete crack surfaces, and thus the need to correspondingly modify or locally adapt the mesh.
Continuous models are typically based on Continuum Damage Mechanics \cite{Krajcinovic:96}, where internal damage variables represent the degradation of material properties due to microcrack formation and growth.
Continuous models have been widely used to model concrete fractures, for instance in \cite{Rots:91,Oliver.Cervera.ea:90} using standard Finite Element Methods (FEM),
in \cite{Ortiz.Leroy.ea:87,Belytschko.Fish.ea:88,Simo.Oliver.ea:93} using FEM with embedded discontinuities,
or in \cite{Wu.Li.ea:15} using extended FEM.
A well-known limitation of these models is their strong sensitivity to the mesh orientation, which can lead to non-physical results.

The starting point of the present paper are \emph{phase-field fracture/damage models}, a particular class of continuous models introduced in \cite{Francfort.Marigo:98,Bourdin.Francfort.ea:00}.
These models introduce a scalar crack phase-field to represent the state of the material, whose evolution is governed by a partial differential equation rather than being postulated a priori.
Phase-field models enable the simulation of complex crack patterns and provide a mathematically consistent framework to account for crack nucleation, propagation, and branching without the need for explicit crack-surface tracking.
In \cite{Francfort.Marigo:98, Bourdin.Francfort.ea:08}, based on Griffith's theory~\cite{Griffith:21}, fracture propagation was conceived as a competition between crack surface energy and bulk elastic energy stored in the material. Later contributions~\cite{Miehe.Welschinger.ea:10,Miehe.Hofacker.ea:10} introduced a more mechanically oriented and thermodynamically consistent formulation, additionally accounting for the irreversibility of crack evolution.
So far, most numerical studies of phase-field models for fracture mechanics have relied on the Finite Element Method (FEM); a non-exhaustive list includes~\cite{Miehe.Hofacker.ea:10,Miehe.Welschinger.ea:10,Heister.Wheeler.ea:15,Ambati.Gerasimov.ea:15,Peerlings.Borst.ea:96,Verhoosel.Borst:13}.  In the most complex situations, however, FEM displays critical limitations in terms of supported meshes, and extended versions can lead to poorly-conditioned algebraic systems.

\smallskip

We devise here what is, to the best of our knowledge, the first HHO discretization of phase-field models.
Originally introduced in~\cite{Di-Pietro.Ern.ea:14,Di-Pietro.Ern:15}, HHO methods can be regarded as a generalization of the classical Crouzeix--Raviart finite element to general polygonal/polyhedral meshes and arbitrary order~\cite{Di-Pietro.Droniou:25}; see~\cite{Di-Pietro.Droniou:20,Cicuttin.Ern.ea:21} for a comprehensive introduction and also \cite{De-Dios.Lipnikov.ea:16} for a related method.
Their distinctive advantage is the support of general polygonal/polyhedral elements, which makes it possible, e.g., %
to treat in a seamless way transition elements in standard meshes, %
to use non-conforming local mesh refinement to increase resolution and %
mesh agglomeration~\cite{Bassi.Botti.ea:12} to reduce the computational cost.
Being constructed from spaces of broken polynomial polynomials on the mesh and its skeleton, HHO methods can also naturally accommodate large variations in the solution.
All of the above features make them particularly suited for fracture modelling.
The method proposed in the present work uses as a starting point:
for the mechanical equilibrium, the classical HHO discretizations of linear elasticity~\cite{Di-Pietro.Ern:15}, modified to incorporate a damage-dependent weight function;
for the evolution of the crack phase-field, the HHO method for scalar diffusion-advection-reaction proposed in~\cite{Di-Pietro.Droniou.ea:15}, adapted to incorporate a history term.
Since we cannot expect the exact solution to be regular, we focus on low polynomial orders, which additionally results in smaller and more manageable algebraic problems.
Time stepping relies on a staggered algorithm, where the mechanical equilibrium and damage-field evolution subproblems are solved independently, combined with suitable strategies for the history field update.
For both subproblems, we additionally use static condensation to further reduce the numerical costs.

Notice that this is not the first time that polyhedral methods have been applied to fracture modelling.
In \cite{Aldakheel.Hudobivnik.ea:18}, the authors propose a Virtual Element discretization of a phase-field model analogous to the one considered here; see \cite{Beirao-da-Veiga.Brezzi.ea:14,Beirao-Da-Veiga.Brezzi.ea:23} for an introduction to the Virtual Elements.
While this method shares several features with ours, it hinges on an entirely different set of degrees of freedom, which resemble (but don't coincide with) those of conforming Lagrange finite elements, and is first-order accurate for both variables;
we also refer to the more recent contributions~\cite{Choi.Chi.ea:23,Liu.Aldakheel.ea:23,Chrysikou.Triantafyllou:25,Leng.Svolos.ea:25,Sharma.Himanshu.ea:25}, all focusing on nodal Virtual Elements.
From the engineering point of view, a potential advantage of our method is that, where no damage is present, it yields equilibrated tractions~\cite{Di-Pietro.Ern:15*1}, and thus more physical solutions.
Another potential advantage is linked to the fact that the pattern of the resulting algebraic problems is smaller, since only elements sharing a face (as opposed to a vertex) are connected.
\smallskip

The rest of the paper is organized as follows.
In Section~\ref{sec:continuous.setting} we briefly recall the phase-field model along with various definitions of the history field.
The HHO scheme is defined in Section~\ref{sec:hho}, while the staggered time stepping scheme used for the solution of the corresponding algebraic problem makes the object of Section~\ref{sec:algorithm}.
Finally, a comprehensive numerical validation is provided in Section~\ref{sec:numerical.tests}, including a comparison with a more standard FEM scheme.


\section{Continuous setting}\label{sec:continuous.setting}

Let $\Omega \subset \Real^3$ be a bounded domain corresponding to a medium in which a fracture propagates.
We denote by $\partial\Omega$ its boundary with outward unit normal vector $\vec n$.
The boundary is decomposed as $\partial\Omega = \partial\Omega_\mathrm{D} \sqcup \partial\Omega_\mathrm{N}$, where $\partial\Omega_\mathrm{D}$ and $\partial\Omega_\mathrm{N}$ are two non-overlapping portions where displacement (Dirichlet) or normal stress (Neumann) conditions are respectively enforced.
We assume small deformations, so that the infinitesimal strain tensor corresponding to a displacement field $\boldsymbol{u} : \Omega \to \Real^3$ is $\tens\varepsilon = \tens\varepsilon(\vec u) = \GRADs \vec{u}$, with $\GRADs$ denoting the symmmetric part of the gradient applied to vector-valued fields.

Throughout the rest of the paper, when dealing with a function $\psi$ of time and space, at any time $t$ for which $\psi$ is defined, we let $\psi(t) \coloneqq \psi(t,\cdot)$ denote the function of space only obtained fixing the first argument at $t$.
We also use the shortcut notation $\dot{\psi}$ to denote the time derivative of $\psi$.

Denoting by $\tF$ the final time, the fracture propagation problem with phase-field approach is formulated as follows:
Given a prescribed boundary displacement $\vec u_\mathrm{D}:(0,\tF] \times \partial\Omega_{\rm D} \to \Real^3$, an initial history field $\calH_0:\Omega\to\Real$, and an initial crack phase-field $\phi_0:\Omega \to [0,1]$ satisfying $\GRAD \phi_0 \cdot \vec n = 0$ on $\partial\Omega$, find the displacement field $\vec u:(0,\tF] \times \Omega \to \Real^3$ and the crack phase-field $\phi:[0,\tF] \times \Omega \to \Real$ such that, for $t \in (0,\tF]$,
\begin{subequations}\label{eq:fracture:strong}
  \begin{alignat}{2}
    -\GRAD \cdot \tens\sigma(\tens\varepsilon(t),\phi(t)) &= \vec 0 \quad &&\text{in }\Omega, \label{eq:fracture:forcebalance} \\
    -\Delta \phi(t) + \frac{\phi(t)}{\ell^2} + \frac{\eta}{G_c\ell}\dot{\phi}(t)
    &= \frac{2(1-\phi(t))}{\ell G_c}\,\calH(t) \quad &&\text{in }\Omega, \label{eq:fracture:phasefield} \\
    \vec u(t) &= \vec u_\mathrm{D}(t) \quad &&\text{on }\partial\Omega_\mathrm{D}, \\
    \tens\sigma(\tens\varepsilon(t),\phi(t)) \cdot \vec n &= \vec 0 \quad &&\text{on }\partial\Omega_\mathrm{N}, \\
    \GRAD \phi(t) \cdot \vec n &= 0 \quad &&\text{on }\partial\Omega,
  \end{alignat}
completed with the initial conditions
\begin{equation}
  \phi(0,\cdot) = \phi_0,\qquad
  \calH(0) = \calH_0.
\end{equation}
\end{subequations}

The stress tensor $\tens\sigma(\tens\varepsilon,\phi)$ is defined by degrading the elastic stress according to the crack phase-field variable, in order to reflect the fact that the stress vanishes in the fractured region:
\[
\tens\sigma(\tens\varepsilon,\phi) \coloneqq g(\phi)\,\tens\sigma_0(\tens\varepsilon).
\]
Above, denoting by $\lambda$ and $\mu$ the Lamé coefficients, $\tens\sigma_0(\tens\varepsilon)$ denotes the standard linear isotropic elasticity stress tensor expressed in terms of the elastic energy $\psi_0$ as
\begin{equation}\label{eq:elastic-energy-stress}
  \text{%
    $\tens\sigma_0(\tens\varepsilon) \coloneqq \frac{\partial\psi_0(\tens\varepsilon)}{\partial\tens\varepsilon} = 2\mu\,\tens\varepsilon + \lambda\,\mathrm{tr}(\tens\varepsilon)\,\tens I$
    where
    $\psi_0(\tens\varepsilon) \coloneqq \frac{\lambda}{2}\big[\mathrm{tr}(\tens\varepsilon)\big]^2 + \mu\,\tens\varepsilon:\tens\varepsilon$,
  }
\end{equation}
with $g$ denoting a degradation function satisfying $g(0)=1$ and $g(1)=0$.
We adopt the common choice $g(\phi) = (1-\phi)^2$.
In \eqref{eq:fracture:phasefield}, $G_c$ is the critical energy release rate from Griffith’s theory and $\ell$ is a length-scale parameter controlling the width of the diffused crack zone.
The term $\eta\,\dot{\phi}$ introduces a viscous regularization with no particular physical meaning, and can be omitted if rate-independent evolution is desired by setting $\eta=0$.

\begin{remark}[The role of time]
  The system~\eqref{eq:fracture:strong} is formulated under quasi-static assumptions.  
  In this setting, the variable $t$ does not represent a physical time, but rather a measure of the progression of the applied load.  
  For this reason, it is often referred to as a \emph{pseudo-time}, and a pseudo-time or time increment simply corresponds to a loading step.
\end{remark}

The first equation \eqref{eq:fracture:forcebalance} expresses a macroscopic equilibrium condition.
Equation~\eqref{eq:fracture:phasefield} is the crack phase-field evolution equation, a micro-balance governing the fracture evolution  in which the right-hand side is the driving force.
The fracture irreversibility is enforced through the \emph{history field} $\calH$, defined as
\begin{equation}\label{eq:historyfield}
  \calH(t) \coloneqq \max_{\tau \in [0,t]} \psi_0(\tens\varepsilon(\tau)),
\end{equation}
which ensures, in particular, that $\calH$ is non-decreasing in time.
Using this definition of $\calH$ leads to the so-called \emph{isotropic formulation}, where both tensile and compressive components of the strain energy equally contribute to fracture growth.

This model may, however, yield unphysical crack propagation under compression due to interpenetration.
To better capture realistic fracture behaviors, it is customary to distinguish tensile (denoted with the superscript ``$+$'') and compressive (denoted with the superscript ``$-$'') contributions to the elastic energy:
\begin{equation}\label{eq:energydecomposition}
  \psi_0(\tens\varepsilon) = \psi_0^+(\tens\varepsilon) + \psi_0^-(\tens\varepsilon),
\end{equation}
and to use only the tensile part in the fracture driving force,
\begin{equation}\label{eq:historyfield:hybrid}
  \calH(t) = \max_{\tau \in [0,t]} \psi_0^+(\tens\varepsilon(\tau)).
\end{equation}
This choice defines the \emph{hybrid formulation}.
Different criteria exist for performing the split \eqref{eq:energydecomposition}.
One possibility consists in using the following spectral decomposition of the strain tensor:
\[
\tens\varepsilon=\sum_{a=1}^3\varepsilon_a\,\vec \normal_a\otimes\vec \normal_a,
\]
where $\varepsilon_a$ ($a=1,2,3$) are the principal strains and $\vec \normal_a$ the corresponding principal directions. Defining
\[
  \text{%
    $\tens\varepsilon^+ \coloneqq \sum_{a=1}^3\varepsilon_a
    ^{\oplus}\,\vec \normal_a\otimes\vec \normal_a$
    and
    $\tens\varepsilon^- \coloneqq \sum_{a=1}^3\varepsilon_a
    ^{\ominus}\,\vec \normal_a\otimes\vec \normal_a$,
  }
\]
with $\alpha^\oplus \coloneqq \frac{|\alpha| + \alpha}{2}$ and $\alpha^\ominus \coloneqq \frac{|\alpha| - \alpha}{2}$ respectively denoting the positive and negative parts of $\alpha$,
the \emph{spectral decomposition} (SP) of the elastic energy reads \cite{Miehe.Welschinger.ea:10}:
\begin{equation}\label{eq:miehe-decomposition}
\text{%
$\psi_0^+(\tens\varepsilon)\coloneq\frac{\lambda}{2} \left( [\tr\tens\varepsilon]^\oplus\right)^2+\mu\,\tens\varepsilon^+:\tens\varepsilon^+$
and $\psi_0^-(\tens\varepsilon)\coloneq\frac{\lambda}{2}\left( [\tr\tens\varepsilon]^\ominus\right)^2+\mu\,\tens\varepsilon^-:\tens\varepsilon^-$.
}
\end{equation}
Notice that both contributions are non-negative.
The evaluation of this energy requires computing the eigenstructure of the strain tensor, which makes the procedure more computationally demanding.
An alternative is the \emph{volumetric–deviatoric decomposition} (VD) of~\cite{Amor.Marigo.ea:09}:
\begin{equation}\label{eq:amor-decomposition}
    \psi_0^+(\tens\varepsilon)\coloneq\frac{1}{2}K \left([\tr\tens\varepsilon]^\oplus\right)^2+\mu\,(\tens\varepsilon':\tens\varepsilon'), \qquad
    \psi_0^-(\tens\varepsilon)\coloneq\frac{1}{2}K \left( [\tr\tens\varepsilon]^\ominus\right)^2,
\end{equation}
where $\tens\varepsilon' \coloneqq \tens\varepsilon - \tfrac{1}{3}\tr(\tens\varepsilon)\tens I$ is the deviatoric component of the strain tensor, and $K = \lambda + \tfrac{2}{3}\mu$ is the bulk modulus.

Problem~\eqref{eq:fracture:strong} with the isotropic definition of the history field \eqref{eq:historyfield} is referred to as the \emph{isotropic formulation}, as the crack is equally driven by both energy components, whereas the version using \eqref{eq:historyfield:hybrid} is referred to as the \emph{hybrid formulation}, which better distinguishes fracture behavior under tension and compression.
The term \emph{hybrid} reflects the intermediate nature of this model, which lies between the isotropic formulation and the fully \emph{anisotropic} one, where the tension–compression split is also applied to the equilibrium equation, a case not considered here.

\begin{remark}[Two-dimensional case]\label{rem:2d}
  The model described above is formulated in three dimensions.
  In the numerical experiments of Section~\ref{sec:numerical.tests}, we consider the two-dimensional case obtained under the plane strain assumption, namely $\varepsilon_{13} = \varepsilon_{23} = \varepsilon_{33} = 0$.
  In what follows, we stick to the three-dimensional terminology for the sake of simplicity.
\end{remark}


\section{HHO Scheme}\label{sec:hho}

We consider a computational mesh $\Mh = \Th \cup \Fh$ covering the domain $\Omega$.
The mesh is made of non-overlapping elements $T \in \Th$, each having a polyhedral shape, and of (planar) polygonal faces $F \in \Fh$, which can be either internal or located on the boundary.
Additional details on the assumptions and properties of such general meshes can be found in \cite[Chapter~1]{Di-Pietro.Droniou:20}.

\subsection{Discrete spaces}\label{sec:hho:discretespaces}

The vector-valued HHO space used for the displacement field is spanned by vector of local functions attached to cells and faces:
\begin{multline*}
  \vVh \coloneq
  \Big\{
  \underline{\vec v}_h = ((\vec v_T)_{T\in\Th}, (\vec v_F)_{F\in\Fh})
  \st
  \\
  \text{$\vec v_T \in \Poly{1}(T)^3$ for all $T\in\Th$ and
        $\vec v_F \in \Poly{1}(F)^3$ for all $F\in\Fh$}
  \Big\},
\end{multline*}
where, for $X \in \Mh$ mesh cell or face, $\Poly{1}(X)$ is spanned by the restriction to $X$ of linear polynomials in the space variables.
Notice that, contrary to finite elements, we do not define here a global space of functions over $\Omega$, but rather work directly on the values of the degrees of freedom, which correspond to polynomial moments of degree up to one over each $X \in \Mh$.
Similarly, the scalar HHO space used for the space discretization of the crack phase-field is defined as
\begin{multline*}
  \Vh\coloneqq
  \Big\{
  \underline{\phi}_h = ((\phi_T)_{T\in\Th}, (\phi_F)_{F\in\Fh})
  \st
  \\
  \text{$\phi_T\in\Poly{0}(T)$ for all $T\in\Th$ and
    $\phi_F\in\Poly{0}(F)$ for all $F\in\Fh$}
  \Big\},
\end{multline*}
with $\Poly{0}(X)$ spanned by constant functions over $X \in \Mh$.
For future use, for any $\underline\phi_h\in\Vh$, we denote by $\phi_h$ the broken polynomial function on $\Omega$ obtained patching cell values, i.e.,
\begin{equation}\label{eq:phi.h}
  \text{
    $(\phi_h)_{|T} \coloneqq \phi_T$ for all $T \in \Th$.
  }
\end{equation}
\begin{remark}[Choice of the polynomial degrees]\label{rem:polynomial.degrees}
  Given the nature of problem~\eqref{eq:fracture:strong}, one cannot expect, in general, high regularity for the exact solution.
  This advocates for the use of low polynomial orders for both the displacement field and the crack phase-field.
  Moreover, in the context of nonlinear problems, low-order schemes often translate into algebraic problems that are easier to solve numerically.
  Notice that the space for displacements is taken of degree $1$ instead of $0$ as this is required for the stability of the standard HHO formulation for elasticity~\cite{Di-Pietro.Ern:15}; see also~\cite{Botti.Di-Pietro.ea:17,Di-Pietro.Droniou:20} on this subject.
  Using the space of degree $0$ would require specific modifications to recover stability. One possibility is to add additional stabilization terms, which results in a method with a larger stencil~\cite{Botti.Di-Pietro.ea:19}. Alternatively, stability can also be achieved without enlarging the stencil by defining the HHO space on suitable subfaces, provided that unisolvency for rigid-body motions is preserved~\cite{Coatleven:23}, though this comes at the price of a larger number of degrees of freedom.
\end{remark}
In what follows, given a cell $T \in \Th$, we will respectively denote by $\vVT$ and $\VT$ the restrictions of $\vVh$ and $\Vh$ to $T$, obtained collecting the components associated with $T$ and its faces.
Hence, for example, $\underline{\phi}_T \in \VT$ denotes the vector $(\phi_T, (\phi_F)_{F \in \FT})$, with $\FT \subset \Fh$ collecting the faces of $T$.

\subsection{Discretization of the mechanical equilibrium equation}\label{sec:hho:mechanics}

The mechanical equilibrium equation~\eqref{eq:fracture:forcebalance} features a linear elastic response weighted by a degradation function depending on the crack phase-field~$\phi$.
The starting point for our formulation is therefore the standard HHO scheme for elasticity \cite{Di-Pietro.Ern:15}.


\subsubsection{Strain reconstruction operator}
To discretize the stress tensor, we introduce an affine strain reconstruction designed by mimicking discrete integration-by-parts formulas.
Specifically, let $T \in \Th$ be a mesh element and let $\Poly{1}(T;\Real_\mathrm{Sym}^{3 \times 3})$ denote the space of symmetric matrix-valued fields with linear polynomial components, i.e.,
\[
\Poly{1}(T;\Real_\mathrm{Sym}^{3 \times 3})
\coloneqq
\left\{
\tens{\tau} =
\begin{pmatrix}
  a & b & c \\
  b & d & e \\
  c & e & f
\end{pmatrix}
\st
a,b,c,d,e,f \in \Poly{1}(T)
\right\}.
\]
Given $\uvec v_T = (\vec v_T, (\vec v_F)_{F \in \FT}) \in \vVT$, the strain reconstruction $\tET{1} \uvec{v}_T \in \Poly{1}(T;\Real_\mathrm{Sym}^{3 \times 3})$ is defined mimicking an integration by parts formula:
\[
\int_T \tET{1} \uvec{v}_T : \tens{\tau}
= - \int_T \vec{v}_T \cdot \DIV \tens{\tau}
+ \sum_{F \in \Fh[T]} \int_F \vec{v}_F \cdot (\tens{\tau} \normal_{TF})
\qquad \text{for all } \tens{\tau} \in \Poly{1}(T; \Real_\mathrm{Sym}^{3 \times 3}),
\]
where $\normal_{TF}$ denotes the normal vector to $F$ pointing out of $T$.
Crucially, this strain reconstruction is one order higher than the one obtained by applying the symmetric gradient operator to the element component $\boldsymbol{v}_T$, making the resulting scheme comparable to a P2 FEM in terms of accuracy.
The associated divergence reconstruction operator is then naturally defined setting
\[
\DT{1} \uvec{v}_T \coloneqq \tr\left(\tET{1} \uvec{v}_T\right).
\]

\subsubsection{Global discrete function $a_h$}

We are now ready to define the global discrete function
$a_h : (\vVh \times \Vh) \times \vVh \to \Real$ such that,
for all $\uvec{u}_h, \uvec{v}_h \in \vVh$ and all $\underline{\phi}_h \in \Vh$,
\begin{equation*}
  a_h((\uvec{u}_h, \underline{\phi}_h), \uvec{v}_h)
  \coloneqq
  \sum_{T\in\Th} g(\phi_T)\, a_T(\uT,\vT),
\end{equation*}
with, for all $T \in \Th$,
\begin{equation}\label{eq:hho:mechanics:aT}
  a_T(\uT,\vT)
  \coloneqq
  \int_T \Big(
    2\mu\,\tET{1}\uvec{u}_T : \tET{1}\uvec{v}_T
    + \lambda\,\DT{1}\uvec{u}_T\,\DT{1}\uvec{v}_T
  \Big)
  + 2\mu\,s_T(\uvec{u}_T, \uvec{v}_T).
\end{equation}
Notice that, when $\phi_T = 0$ for all $T \in \Th$ (unfractured medium), $a_h$ reduces to the standard HHO bilinear form for linear elasticity.
The local bilinear form~\eqref{eq:hho:mechanics:aT} is composed of two contributions, respectively in charge of consistency and stability.
The stabilization term $s_T$ is designed so that:
\begin{itemize}
  \item[(S1)] $a_T$ vanishes only when one of its arguments corresponds to the interpolation of a rigid-body motion;
  \item[(S2)] $s_T$ vanishes whenever the exact displacement is a quadratic function.
\end{itemize}
In order to fulfill these properties, $s_T$ must have a very specific dependence on its arguments, which requires the introduction of a reconstructed quadratic displacement field.
Denoting by $\Poly{2}(T)$ the space of quadratic polynomials on $T$, for a given $\uvec v_T = (\vec v_T, (\vec v_F)_{F \in \FT}) \in \vVT$, we define $\vpT{2} \uvec{v}_T \in \Poly{2}(T)^3$ such that
its symmetric gradient is the projection of $\tET{1}$ on symmetric gradients of quadratic polynomials
\begin{subequations}\label{eq:elasticity:vpT}
\begin{equation}
  \int_T\big(\GRADs \vpT{2} \uvec{v}_T - \tET{1} \uvec{v}_T\big) : \GRADs \vec{w} = 0
  \qquad \text{for all } \vec{w} \in \Poly{2}(T)^3,
\end{equation}
and we fix rigid-body motions additionally requiring that
\begin{equation}\label{eq:vpT:reconstruction}
  \int_T \vpT{2} \uvec{v}_T = \int_T \vec{v}_T,
  \qquad
  \int_T \GRADss \vpT{2} \uvec{v}_T
  = \frac{1}{2} \sum_{F \in \Fh[T]} \int_F
  \left(
  \vec{v}_F \otimes \normal_{TF}
  - \normal_{TF} \otimes \vec{v}_F
  \right),
\end{equation}
\end{subequations}
with $\GRADss$ denoting the skew-symmetric part of the gradient applied to vector-valued fields.
This quadratic displacement field can be equivalently characterized
as providing the best approximation of $\tET{1} \uvec{v}_T$ in $\GRADs \Poly{2}(T)^3$:
\[
\vpT{2} \uvec{v}_T
= \argmin_{\text{$\vec w \in \Poly{2}(T)^3$ satisfying~\eqref{eq:vpT:reconstruction}}} \frac{1}{2} \norm[L^2(T)^{3\times 3}]{\GRADs \vec{w} - \tET{1} \uvec{v}_T}^2.
\]
It can be proved that properties (S1)--(S2) can only be fulfilled if $s_T$ is a least-squares penalization of the following quantities:
\[
\vec{\delta}_T^1 \uvec{v}_T
\coloneq
\boldsymbol{\pi}_T^1(\vpT{2} \uvec{v}_T - \vec{v}_T),
\qquad
\vec{\delta}_{TF}^1 \uvec{v}_T
\coloneq
\boldsymbol{\pi}_F^1(\vpT{2} \uvec{v}_T - \vec{v}_F)
\quad \text{for all }\, F \in \Fh[T],
\]
where, for $X \in \Mh$, given a smooth-enough displacement field $\vec v$ on $X$, $\vec{\pi}_X^1 \vec{v}$ is the polynomial in $\Poly{1}(X)^3$ providing the best approximation of $\vec v$ in the $L^2$-norm.
Examples of stabilization are
\[
s_T(\uvec u_T, \uvec v_T)
= \sum_{F \in \FT} h_F^{-1} \int_F \left(\vec{\delta}_{TF}^1 \uvec u_T - \vec{\delta}_T^1 \uvec u_T\right) \cdot \left( \vec{\delta}_{TF}^1 \uvec v_T - \vec{\delta}_T^1 \uvec v_T\right)
\] 
or
\[
s_T(\uvec u_T, \uvec v_T)
= h_T^{-2} \int_T \vec{\delta}_T^1 \uvec u_T \cdot \vec{\delta}_T^1 \uvec v_T
+ \sum_{F \in \FT} h_F^{-1} \int_F \vec{\delta}_{TF}^1 \uvec u_T \cdot \vec{\delta}_{TF}^1 \uvec v_T,
\]
where $h_T$ denotes the diameter of $T$ and $h_F$ the diameter of $F$.

\subsection{Discretization of the crack phase-field evolution equation}\label{sec:hho:phasefield}

The crack phase-field equation~\eqref{eq:fracture:phasefield}
can regarded as an evolutive diffusion–reaction equation for the crack phase-field.
Linear standalone equations of this kind have been extensively studied, and the corresponding HHO discretizations are well-established.
We specifically use here as a starting point the method of~\cite{Di-Pietro.Droniou.ea:15} (see also \cite[Chapter~3]{Di-Pietro.Droniou:20}), which is modified to accommodate the specific fracture driving terms arising in our model.

\subsubsection{Local reconstructions}\label{sec:hho:phasefield:reconstructions}

As in the mechanical problem, we introduce a higher-order polynomial reconstruction operator for the crack phase-field.
Specifically, for any cell $T \in \Th$, given $\underline \phi_T = (\phi_T, (\phi_F)_{F \in \FT}) \in \VT$, the affine phase-variable reconstruction $\pT{1} \underline{\phi}_T$, obtained in a similar spirit as the displacement reconstruction~\eqref{eq:elasticity:vpT}, is such that
\[
\text{%
  $\pT{1} \underline{\phi}_T = \phi_T + \GRAD\pT{1} \underline{\phi}_T \cdot  \big(\vec{x} - \bar{\vec{x}}_T\big)$
  with
  $\GRAD \pT{1} \underline{\phi}_T = \frac{1}{|T|} \sum_{F \in \Fh[T]} |F|\, \phi_F\, \normal_{TF},$
}
\]
with $\overline{x}_T$ denoting the center of mass of $T$, $|T|$ its volume, and $|F|$ the area of $F$.

\subsubsection{Global discrete function $b_h$}\label{sec:hho:phasefield:globalproblem}

Given a bounded function $\mathscr{H} : \Omega \to \Real$, we define $b_h(\mathscr{H};\, \cdot,\cdot) : \Vh \times \Vh \to \Real$
such that, for all $\underline{\phi}_h,\underline{\chi}_h \in \Vh$,
\begin{equation}\label{eq:ah:pf}
  b_h(\mathscr{H};\underline{\phi}_h,\underline{\chi}_h)
  \coloneqq \sum_{T \in \Th} 
  \int_T \left[
    \GRAD \pT{1} \underline{\phi}_T \cdot \GRAD \pT{1} \underline{\chi}_T
    + \left(
    \frac{1}{\ell^2} + \frac{2}{\ell G_c}\mathscr{H}
    \right) \phi_T \, \chi_T
    \right]
    + \sum_{T \in \Th} j_T(\underline{\phi}_T, \underline{\chi}_T).
\end{equation}
The first contribution in~\eqref{eq:ah:pf} accounts for the diffusion and reaction contributions, whereas the second one is a stabilizing term in the same spirit as in the mechanical equilibrium problem, obtained setting, for all $T \in \Th$,
\[
j_T(\uphi,\uchi)
\coloneqq \sum_{F\in\Fh[T]}\frac{1}{h_T |F|}\bigg(\int_F (\pT{1} \uphi - \phi_F)\bigg) \bigg(\int_F (\pT{1} \uchi - \chi_F)\bigg).
\]

\subsubsection{Discrete history field variable}

Let $T\in\Th$ be a mesh element.
The time-dependent quantity $\calH_T$ represents the discrete counterpart of the history field~\eqref{eq:historyfield}, i.e., the maximum of the elastic energy density over time.
At a point in time $t$, for the isotropic formulation, we set
\[
\calH_T(t) \coloneqq \max_{\tau \in [0,t]} \psi_0(\tET{1}\uvec u_T(\tau)),
\]
with $\psi_0$ defined by~\eqref{eq:elastic-energy-stress}, while,
for the hybrid formulation, we let
\[
\calH_T(t) \coloneqq \max_{\tau \in [0,t]} \psi^+_0(\tET{1}\uvec{u}_T(\tau)),
\]
where $\psi^+_0$ is given by~\eqref{eq:miehe-decomposition} or~\eqref{eq:amor-decomposition}, depending on the selected decomposition.
At the global level, we then define the function of time $\calH_h : [0,\tF] \times \Omega \to \Real$ such that, for $t \in [0,\tF]$,
\[
(\calH_h(t))_{|T} \coloneqq \calH_T(t) \quad \text{for all }\, T \in \Th.
\]

\subsection{Space semi-discrete problem}\label{sec:hho:discreteproblem}

  We define the following subspaces of $\vVh$ respectively incorporating non-homogeneous and homogeneous boundary conditions  on $\partial\Omega_{\mathrm{D}}$:
\[
\begin{aligned}
  \vVhD
  &\coloneqq
  \left\{
  \underline{\vec v}_h \in \vVh \st
  \text{
    $\vec v_F = \boldsymbol{\pi}_F^1 \vec u_{\mathrm{D}}\quad $for all $F \in \Fh$ such that $F \subset \partial \Omega_{\mathrm{D}}$
  }%
  \right\},
  \\
  \vVhZ
  &\coloneqq
  \left\{
  \underline{\vec v}_h \in \vVh \st
  \text{
    $\vec v_F = \vec 0\quad $for all $F \in \Fh$ such that $F \subset \partial \Omega_{\mathrm{D}}$
  }%
  \right\}.
\end{aligned}
\]
Notice that $\vec u_{\mathrm{D}}$ is a function of time, and so is therefore $\vVhD$.
Hence, $\uvec v_h \in \vVhD$ means that $\uvec v_h$ is a function of time with values in $\vVh$ and such that the boundary condition is satisfied at each time $t \in (0,\tF]$.
The space semi-discrete problem, counterpart of \eqref{eq:fracture:strong}, reads as follows:
Find $\uvec{u}_h \in \vVhD$ and $\underline{\phi}_h:[0,\tF]\to\Vh$ such that, for $t \in (0,\tF]$,
\begin{subequations}\label{eq:hho:fracture:discrete}
  \begin{alignat}{2}\label{eq:hho:fracture:discrete:momentum}
    a_h((\uvec{u}_h(t), \underline{\phi}_h(t)), \uvec{v}_h) &= 0
    && \text{for all }\, \uvec{v}_h \in \vVhZ,
    \\ \label{eq:hho:fracture:discrete:phase-field}
    \int_\Omega \frac{\eta}{\ell G_c}\,\dot{\phi}_h(t)\,\chi_h
    + b_h(\calH_h(t); \underline{\phi}_h(t), \underline{\chi}_h)
    &=\int_\Omega \frac{2}{\ell G_c}\,\calH_h(t)\,\chi_h
    \qquad && \text{for all }\, \underline{\chi}_h \in \Vh,
  \end{alignat}
  where the functions $\phi_h$ and $\chi_h$ in~\eqref{eq:hho:fracture:discrete:phase-field} are defined from $\underline{\phi}_h$ and $\underline{\chi}_h$ according to~\eqref{eq:phi.h},
  completed with the initial conditions
  \[
  \underline{\phi}_h(0) = \left(
  \left(\frac{1}{|T|}\int_T \phi_0 \right)_{T \in \Th},
  \left(\frac{1}{|F|}\int_F \phi_0 \right)_{F \in \Fh}
  \right),\qquad
  \text{
    $\calH_T(0) = \frac{1}{|T|} \int_T \calH_0$ for all $T \in \Th$.
  }
  \]
\end{subequations}


\section{Staggered time stepping scheme}\label{sec:algorithm}

Problem~\eqref{eq:hho:fracture:discrete} is a strongly non-linear system of ordinary differential equations the solution of which requires carefully tailored resolution strategies.  We describe hereafter the one used in the numerical simulations of Section~\ref{sec:numerical.tests}.

\subsection{Time discretization and non-linear iterations}\label{sec:algorithm:basic}

The time discretization hinges on a uniform partition of size $\tau > 0$ of the time interval $[0,\tF]$:
\[
\text{
  $0 = t_0 < t_1 < \ldots < t_N = \tF$ with $t_{n+1} - t_n = \tau$ for $n = 0,\ldots, N-1$.
}
\]
At each time step $t_n$, for $n = 1, \ldots, N$, we let $\vVhDn{n} \coloneqq \vVhD(t^n)$ and denote by $(\uvec u_h^n,\underline{\phi}_h^n,\calH_T^n)$ the approximation of $(\uvec u_h(t_n),\underline{\phi}(t_n),\calH_h(t_n))$.
The non-linear algebraic system to advance from time step $t_n$ to time step $t_{n+1}$ is solved by means of the iterative method described hereafter.
We initialize the crack phase-field and history field using the values at the previous time step:
\[
\underline\phi_h^{n+1,0} = \underline\phi_h^n, \qquad
\calH_h^{n+1,0} = \calH_h^n.
\]
Subsequently, for each non-linear iteration $m=0,1,\ldots$, we compute the updated discrete fields $\uvec u_h^{n+1,m+1}$, $\underline\phi_h^{n+1,m+1}$, and $\calH_h^{n+1,m+1}$ from the fields $\underline\phi_h^{n+1,m}$ and $\calH_h^{n+1,m}$ available from the previous iteration as follows:
\begin{enumerate}
\item \textbf{Mechanical equilibrium problem.} First, we obtain $\uvec{u}_h^{n+1,m+1} \in\vVhDn{n+1}$ from the known crack phase-field $\underline\phi_h^{n+1,m}$ by solving the mechanical equilibrium equation
  \begin{equation}\label{eq:discrete:mechanical}
    a_h((\uvec{u}_h^{n+1,m+1}, \underline{\phi}_h^{n+1,m}), \uvec{v}_h) = 0
    \quad\text{for all }\, \uvec{v}_h \in \vVhZ;
  \end{equation}
\item \textbf{History field update.} Next, for all $T \in \Th$, we update the history field $\calH_T^{n+1,m+1}$ using  $\uvec{u}_T^{n+1,m+1}$, $\uvec{u}_T^{n+1,m}$, and, possibly, $\calH_T^{n+1,m}$.

\item \textbf{Crack phase-field evolution problem.} Finally, we solve the crack phase-field evolution equation to obtain $\underline \phi_h^{n+1,m+1} \in \Vh$, using the updated $\calH_h^{n+1,m+1}$ and a backward Euler scheme:
  \begin{multline}\label{eq:phase-field:discrete}
    \int_\Omega \frac{\eta}{\ell G_c}\,\frac{\phi_h^{n+1,m+1}}{\tau}\,\chi_h
    + b_h(\calH_h^{n+1,m+1}; \underline{\phi}_h^{n+1,m+1}, \underline{\chi}_h)
    \\
    = \int_\Omega \Big(
    \frac{2}{\ell G_c}\,\calH_h^{n+1,m+1}
    + \frac{\eta}{\ell G_c}\,\frac{\phi_h^n}{\tau}
    \Big)
    \,\chi_h
  \quad\text{for all }\, \underline{\chi}_h \in \Vh.
  \end{multline}
\end{enumerate}
We exit the non-linear iterations when the relative increment of both the displacement field and crack phase-field between two consecutive iterations is under a user-defined threshold.
Some remarks are of order.

\begin{remark}[The case of volumetric-deviatoric energy split]
  When using the volumetric–deviatoric energy split \eqref{eq:amor-decomposition} for the tensile and compressive parts of the strain energy, both subproblems~\eqref{eq:discrete:mechanical} and~\eqref{eq:phase-field:discrete} are linear.
\end{remark}

\begin{remark}[Static condensation]
  In both subproblems, the global systems are solved using static condensation: the cell unknowns are locally eliminated, yielding condensed systems that involve only face unknowns. Once the global problem is solved, the cell unknowns are recovered solving an inexpensive equation. This procedure significantly reduces the computational cost.
\end{remark}

\subsection{History field update}\label{sec:H-implementation}

Let a cell $T \in \Th$, a time step $n$, and a non-linear iteration $m$ be given.
We discuss here the history field update in Step 2 of the algorithm presented in Section~\ref{sec:algorithm:basic}.
We assume that we have stored the value $\tens\varepsilon_T^{n+1,m} \in \Poly{1}(T; \Real_{\rm Sym}^{3\times 3})$ of the strain reconstruction corresponding to the largest elastic energy attained so far.

Given a function $f:T\to\Real$, let $f^\mathrm{max}$ denote its maximum over the quadrature nodes of $T$.
The value of the discrete strain tensor associated with the largest energy is then updated as
\[
\tens\varepsilon_T^{n+1,m+1} =
\begin{cases}
  \tET{1}\uvec{u}_T^{n+1,m+1}
  & \text{if } \psi_0^{*,\,\mathrm{max}}(\tET{1}\uvec{u}_T^{n+1,m+1})
  > \calH_T^{n+1,m,\,\mathrm{max}},\\[3mm]
  \tens\varepsilon_T^{n+1,m}
  & \text{otherwise,}
\end{cases}
\]
where $\psi_0^*$  denotes either $\psi_0$ for the isotropic formulation or $\psi_0^+$ for the hybrid one.
Notice that another option consists in replacing $\calH_T^{n+1,m,\,\mathrm{max}}$ with $\calH_T^{n,\,\mathrm{max}}$ in the above equation.
Finally, the cell history field is updated setting
\[
\calH_T^{n+1,m+1} = \psi^*_{0}(\tens\varepsilon_T^{n+1,m+1}).
\]

\begin{remark}[Update at quadrature nodes]
  The above strategy to update the history field yields a polynomial $\calH_T^{n+1,m+1}$ within each cell, which can facilitate the transfer of this quantity between meshes when local mesh refinement is used.
  An alternative consists in first replacing the integrals in~\eqref{eq:phase-field:discrete} with quadratures and then storing and updating the history field at each quadrature node.
  In the numerical experiments we performed, both strategies essentially yields the same results.
\end{remark}


\section{Numerical simulations}\label{sec:numerical.tests}

We tested the proposed scheme on two standard benchmark problems in fracture mechanics.
Both the isotropic and the hybrid formulations were implemented, and for the hybrid one we considered both the spectral \eqref{eq:miehe-decomposition} and the volumetric-deviatoric \eqref{eq:amor-decomposition} energy decompositions.
Our implementation is based on the open source \texttt{HArDCore} library\footnote{\url{https://github.com/jdroniou/HArDCore}}.

\subsection{Traction test}

\begin{figure}
  \centering
  \includegraphics[width=0.8\textwidth]{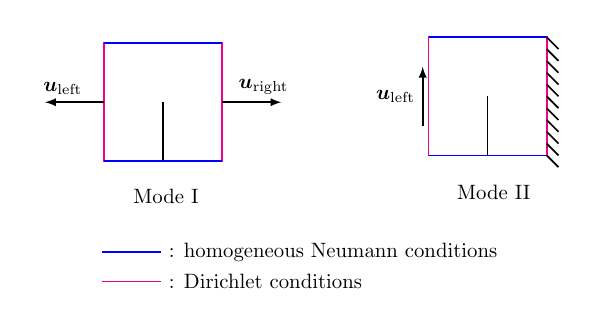}
  \caption{Domain and boundary conditions for the traction and shear test.}
  \label{fig:domain}
\end{figure}

The first benchmark corresponds to the so-called \emph{mode I} tensile test.
We consider a square domain of $\SI{1}{\milli\metre}$ side length containing an initial vertical notch of length $\SI{0.5}{\milli\metre}$, as shown in Figure~\ref{fig:domain}.
A zero displacement is imposed on the right side, while on the left side a monotonically increasing displacement is enforced in the direction normal to the edge, pointing outward.
The prescribed displacement starts from zero and increases linearly with an increment of $\SI{1e-5}{\milli\metre}$ per time step for the first $500$ steps, and then with a smaller increment of $\SI{1e-6}{\milli\metre}$ per time step to better capture the sudden crack propagation.
The domain is discretized using a refined family of uniform conforming triangular meshes, and the phase-field regularization length is set to $\ell = \SI{0.0075}{\milli\metre}$.
The Lamé parameters are $\mu = \SI{80.77}{\kilo\newton\per\milli\metre\squared}$ and $\lambda = \SI{121.15}{\kilo\newton\per\milli\metre\squared}$, while the critical energy release rate is set to $G_c = \SI{2.7e-3}{\kilo\newton\per\milli\metre}$.

\begin{figure}[htbp]
  \centering
  \begin{subfigure}[t]{0.215\textwidth}
    \includegraphics[width=\textwidth]{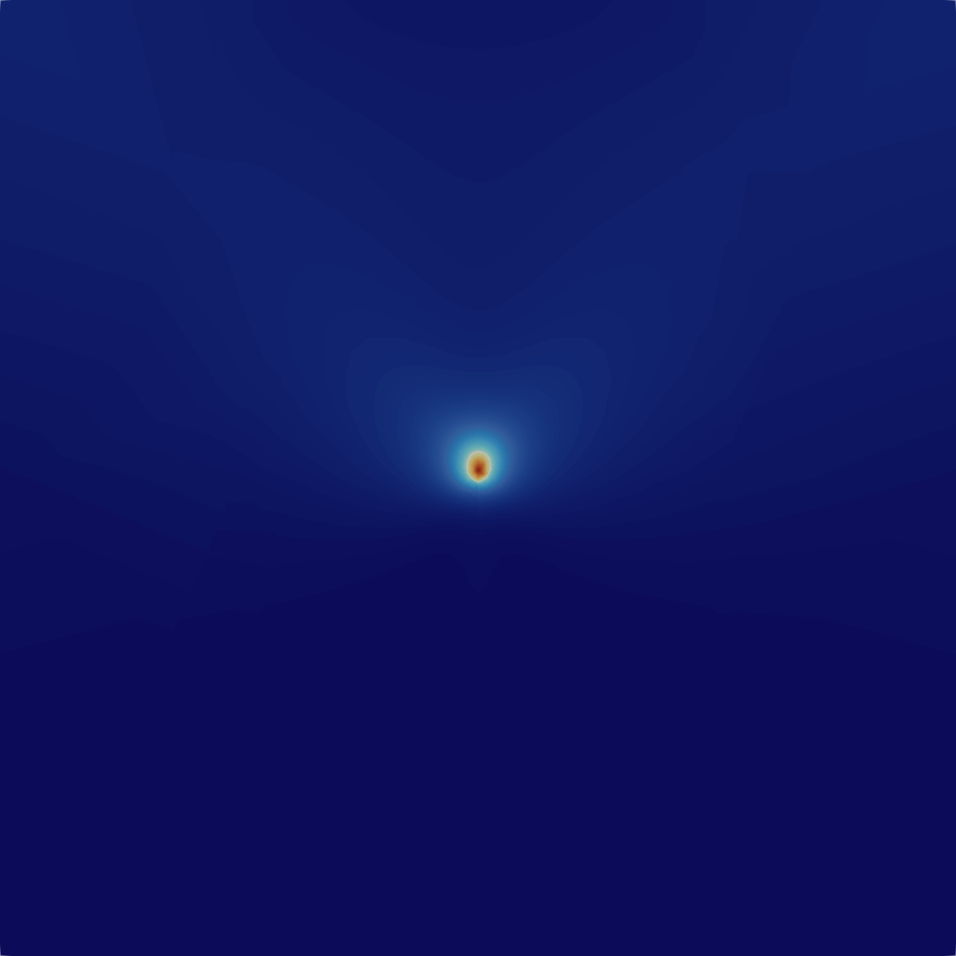}
  \end{subfigure}
  \hspace{0.0001cm}
  \begin{subfigure}[t]{0.215\textwidth}
    \includegraphics[width=\textwidth]{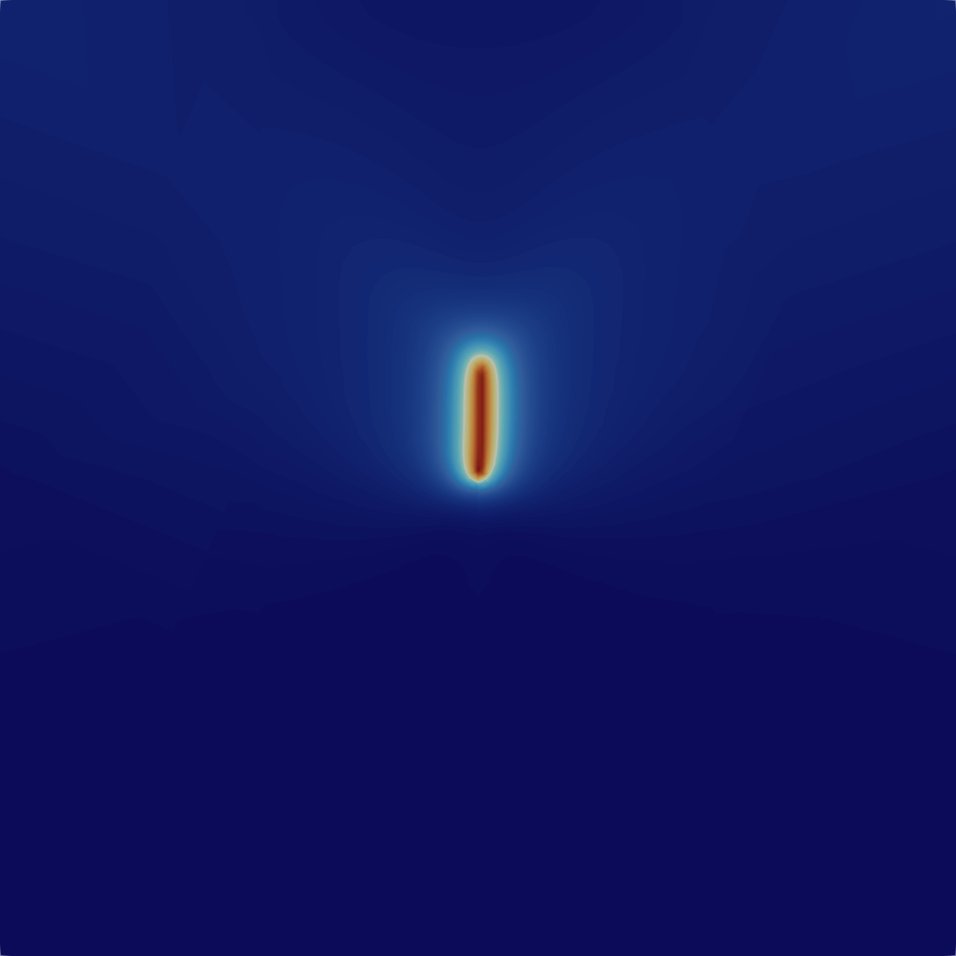}
  \end{subfigure}
  \hspace{0.0001cm}
  \begin{subfigure}[t]{0.215\textwidth}
    \includegraphics[width=\textwidth]{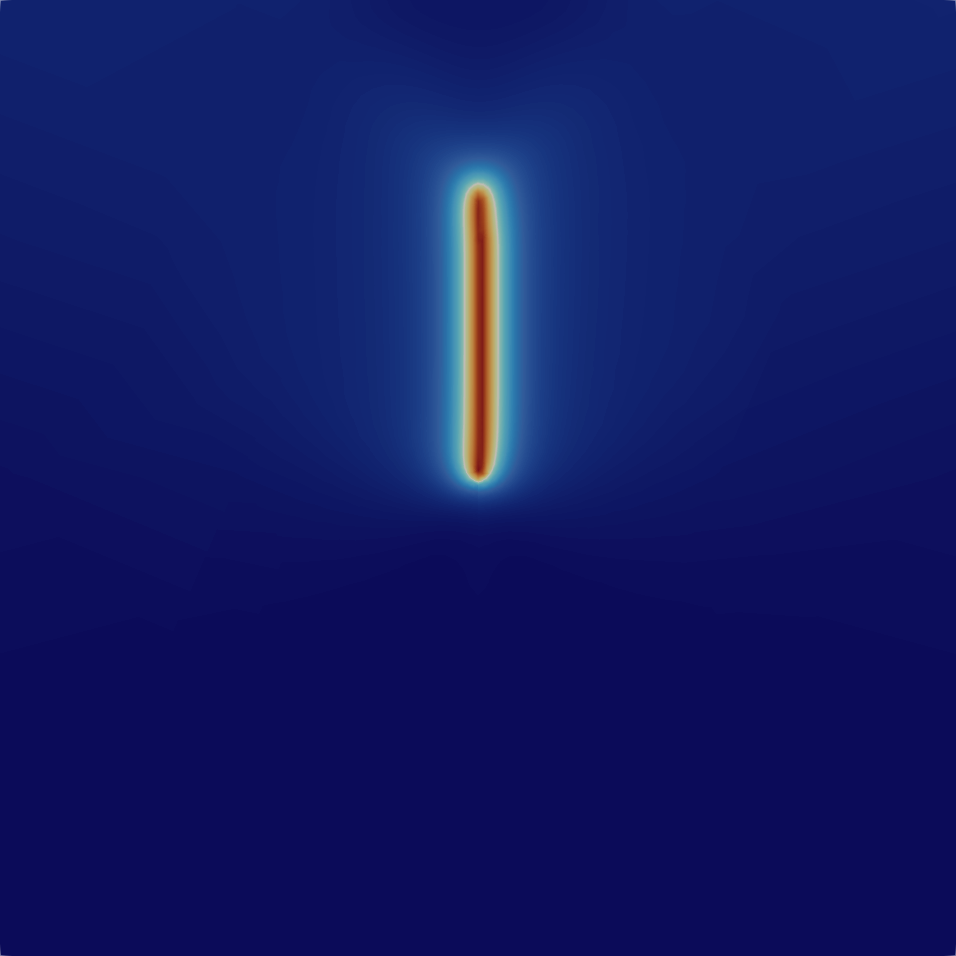}
  \end{subfigure}
  \hspace{0.0001cm}
  \begin{subfigure}[t]{0.215\textwidth}
    \includegraphics[width=\textwidth]{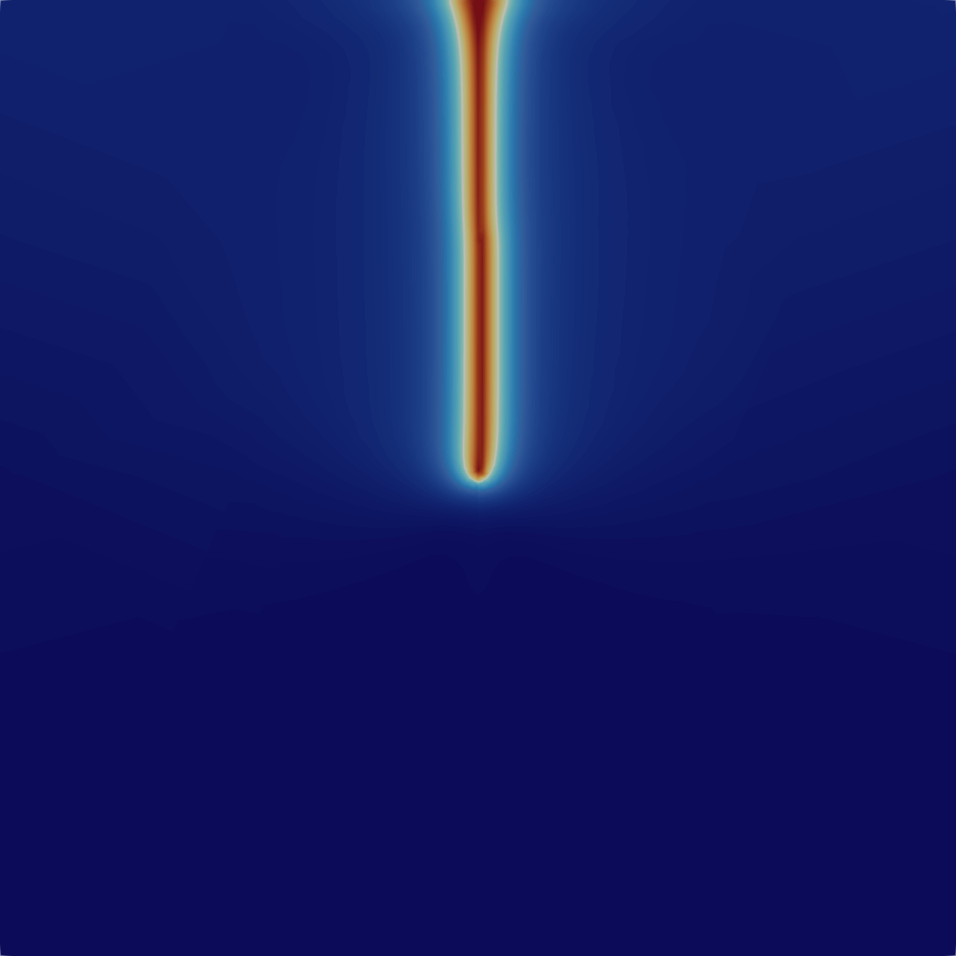}
  \end{subfigure}
  \hspace{0.001cm}
  \begin{subfigure}[t]{0.07\textwidth}
    \vspace{-3.2cm}
    \includegraphics[width=\textwidth]{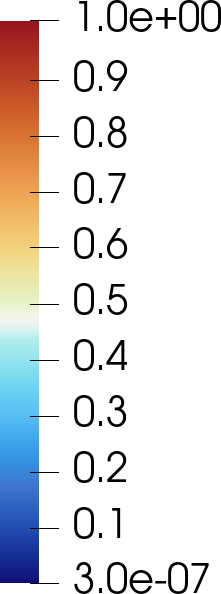}
  \end{subfigure}

  \vspace{0.15cm}
  \begin{subfigure}[t]{0.215\textwidth}
    \includegraphics[width=\textwidth]{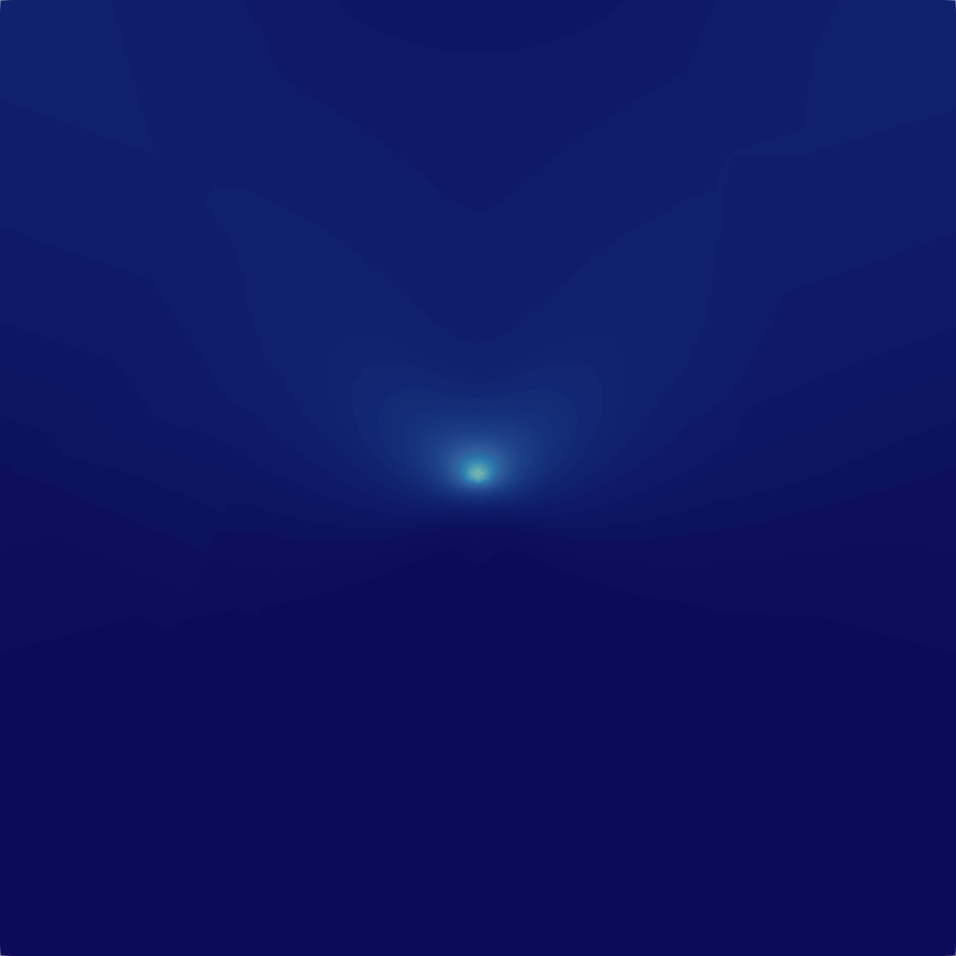}
  \end{subfigure}
  \hspace{0.0001cm}
  \begin{subfigure}[t]{0.215\textwidth}
    \includegraphics[width=\textwidth]{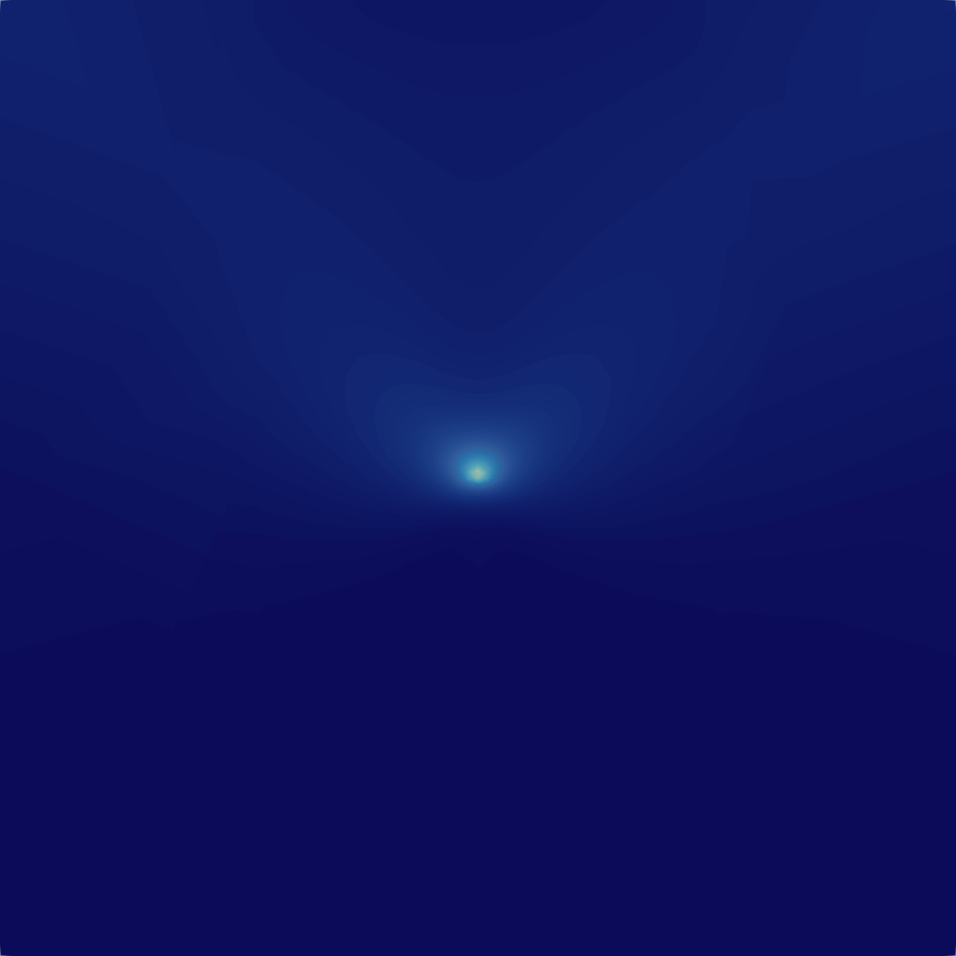}
  \end{subfigure}
  \hspace{0.0001cm}
  \begin{subfigure}[t]{0.215\textwidth}
    \includegraphics[width=\textwidth]{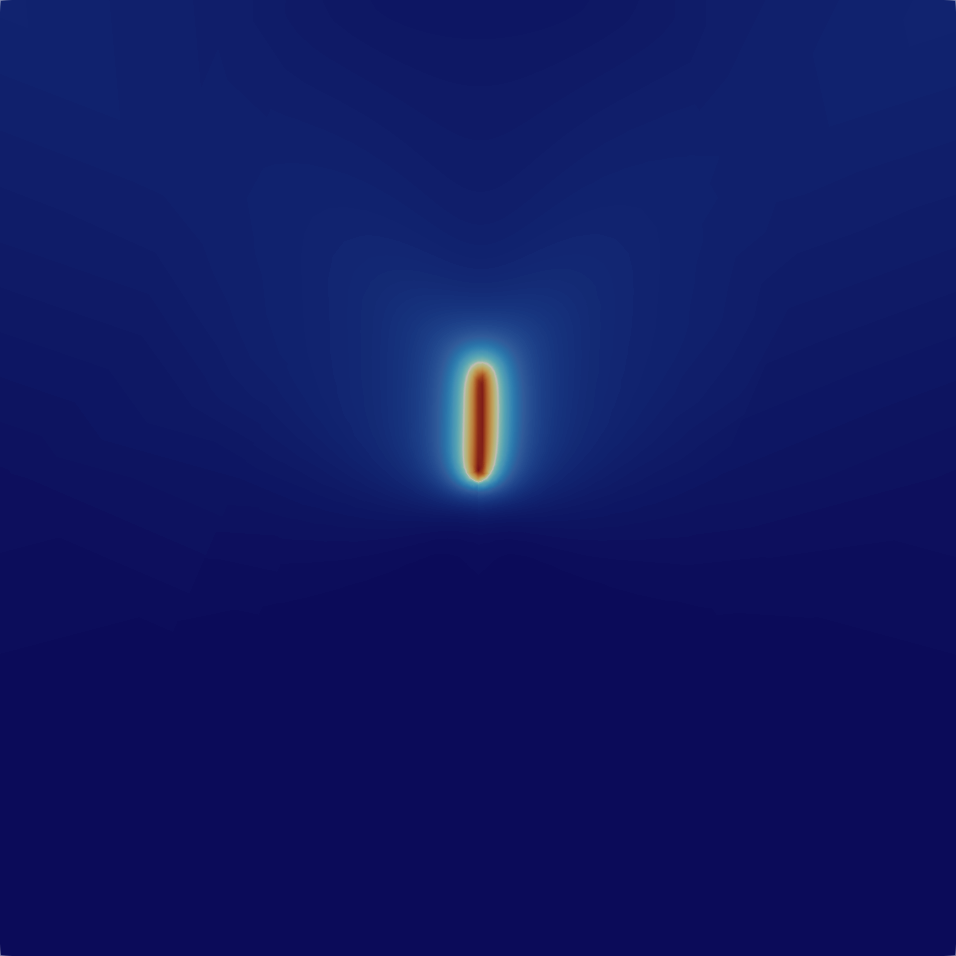}
  \end{subfigure}
  \hspace{0.0001cm}
  \begin{subfigure}[t]{0.215\textwidth}
    \includegraphics[width=\textwidth]{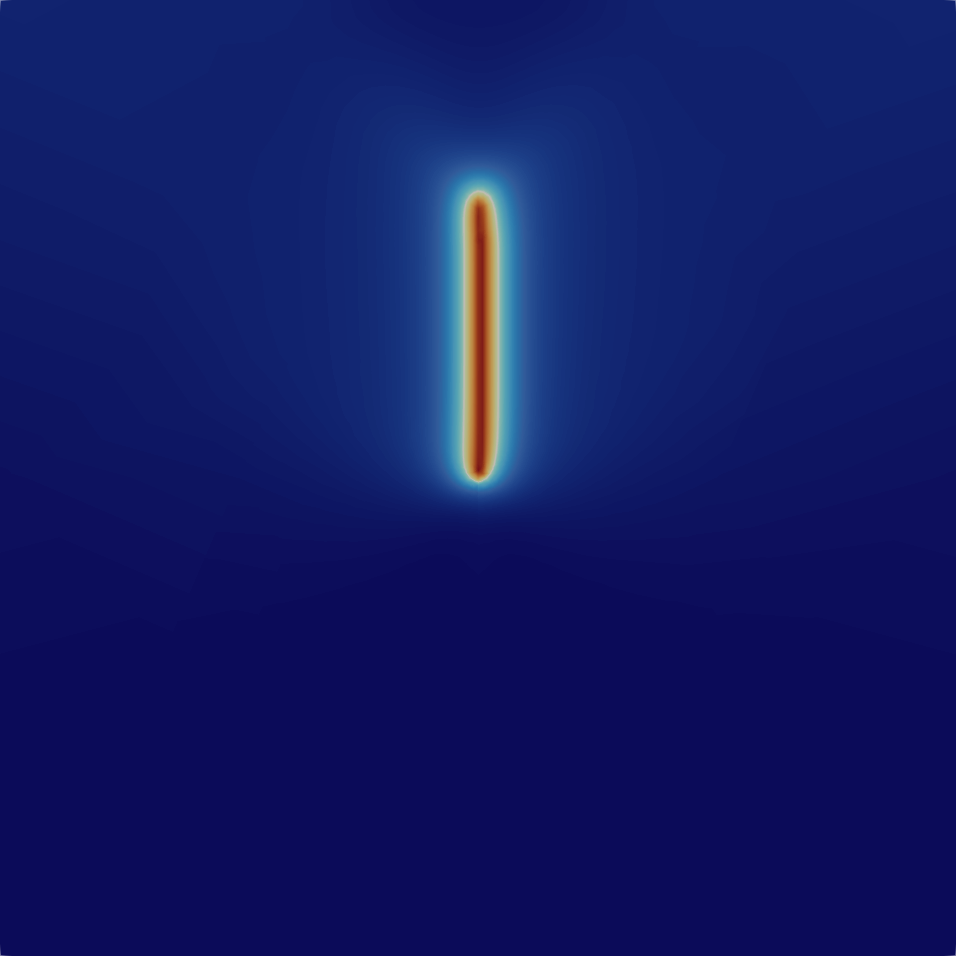}
  \end{subfigure}
  \hspace{0.001cm}
  \begin{subfigure}[t]{0.07\textwidth}
    \vspace{-3.2cm}
    \includegraphics[width=\textwidth]{Figures/hho_simulations/legend.png}
  \end{subfigure}

  \vspace{0.15cm}
  \begin{subfigure}[t]{0.215\textwidth}
    \includegraphics[width=\textwidth]{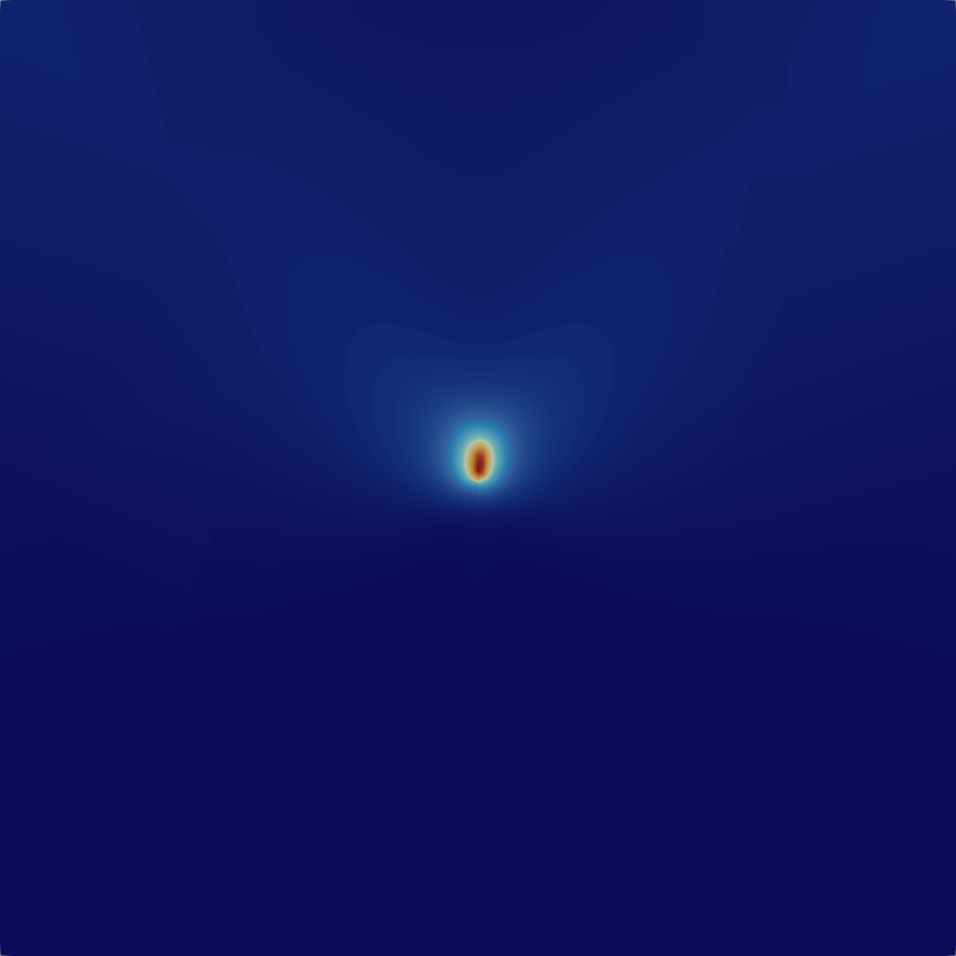}
    \caption*{\scriptsize $t=60$}
  \end{subfigure}
  \hspace{0.0001cm}
  \begin{subfigure}[t]{0.215\textwidth}
    \includegraphics[width=\textwidth]{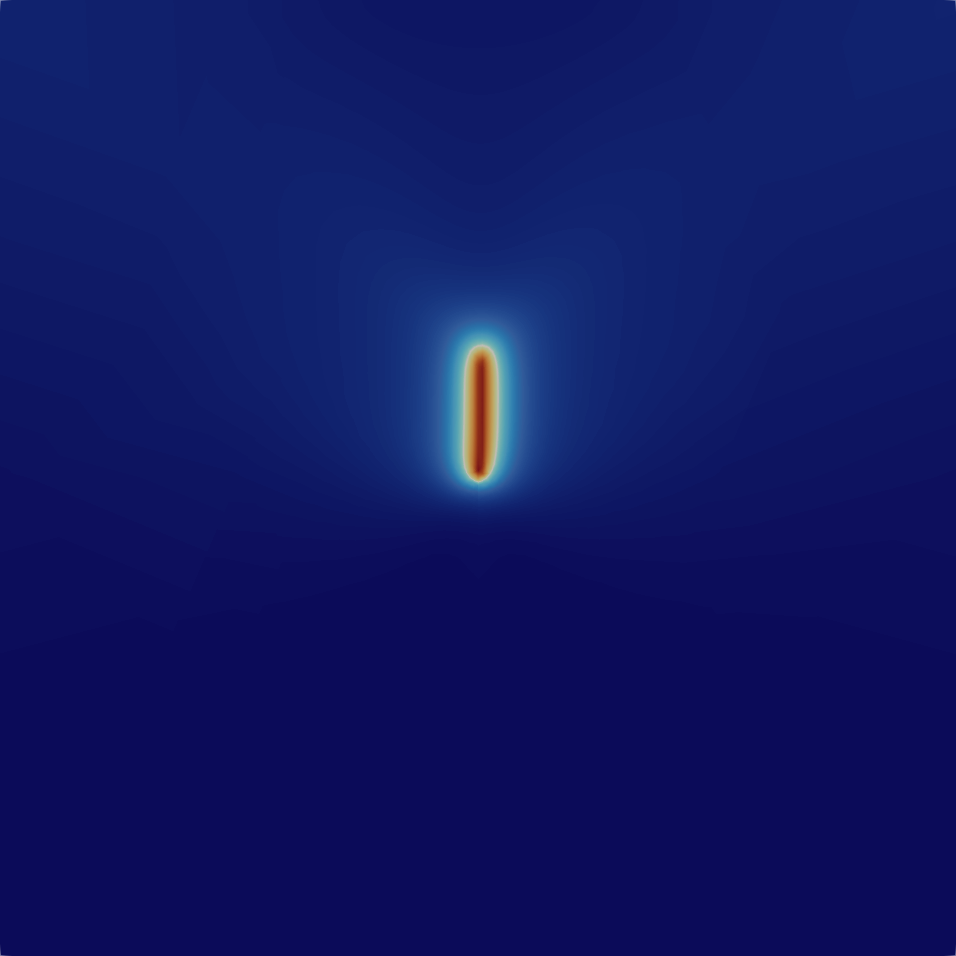}
    \caption*{\scriptsize $t=70$}
  \end{subfigure}
  \hspace{0.0001cm}
  \begin{subfigure}[t]{0.215\textwidth}
    \includegraphics[width=\textwidth]{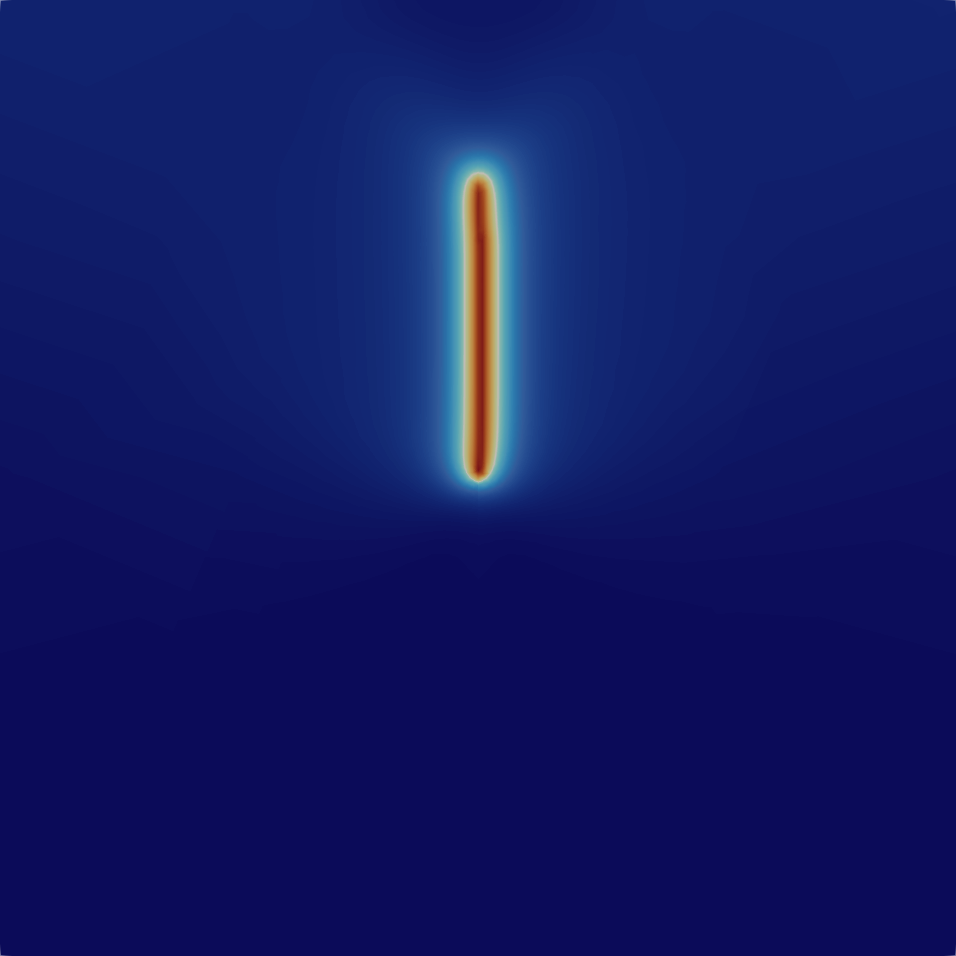}
    \caption*{\scriptsize $t=85$}
  \end{subfigure}
  \hspace{0.0001cm}
  \begin{subfigure}[t]{0.215\textwidth}
    \includegraphics[width=\textwidth]{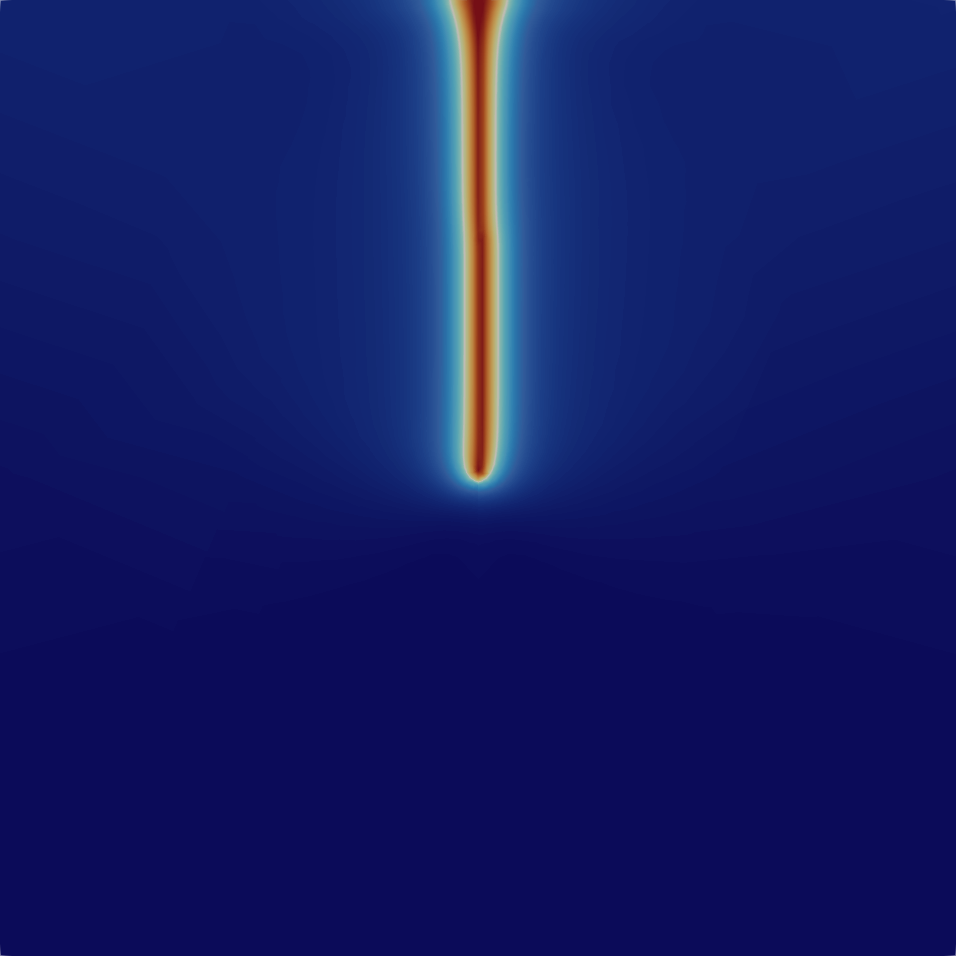}
    \caption*{\scriptsize $t=100$}
  \end{subfigure}
  \hspace{0.001cm}
  \begin{subfigure}[t]{0.07\textwidth}
    \vspace{-3.2cm}
    \includegraphics[width=\textwidth]{Figures/hho_simulations/legend.png}
  \end{subfigure}

  \caption{Evolution of the crack phase-field for different formulations (from top to bottom: isotropic, hybrid–VD, hybrid–SP) at increasing time steps, using a uniform mesh with element size $h = \SI{0.01}{\milli\metre}$}
  \label{fig:phasefield-evolution-modeI}
\end{figure}

The evolution of the crack phase-field $\phi$ at different time steps is shown in Figure~\ref{fig:phasefield-evolution-modeI} for the three formulations implemented with the HHO method: isotropic, hybrid-VD (with energy decomposition~\eqref{eq:amor-decomposition}), and hybrid-SP (with energy decomposition~\eqref{eq:miehe-decomposition}).
The three formulations exhibit very similar behavior in this test, with a crack propagating straight from the initial notch to the top boundary, as expected from physical experiments and the literature.

\begin{figure}
  \centering
  \includegraphics[width=\textwidth]{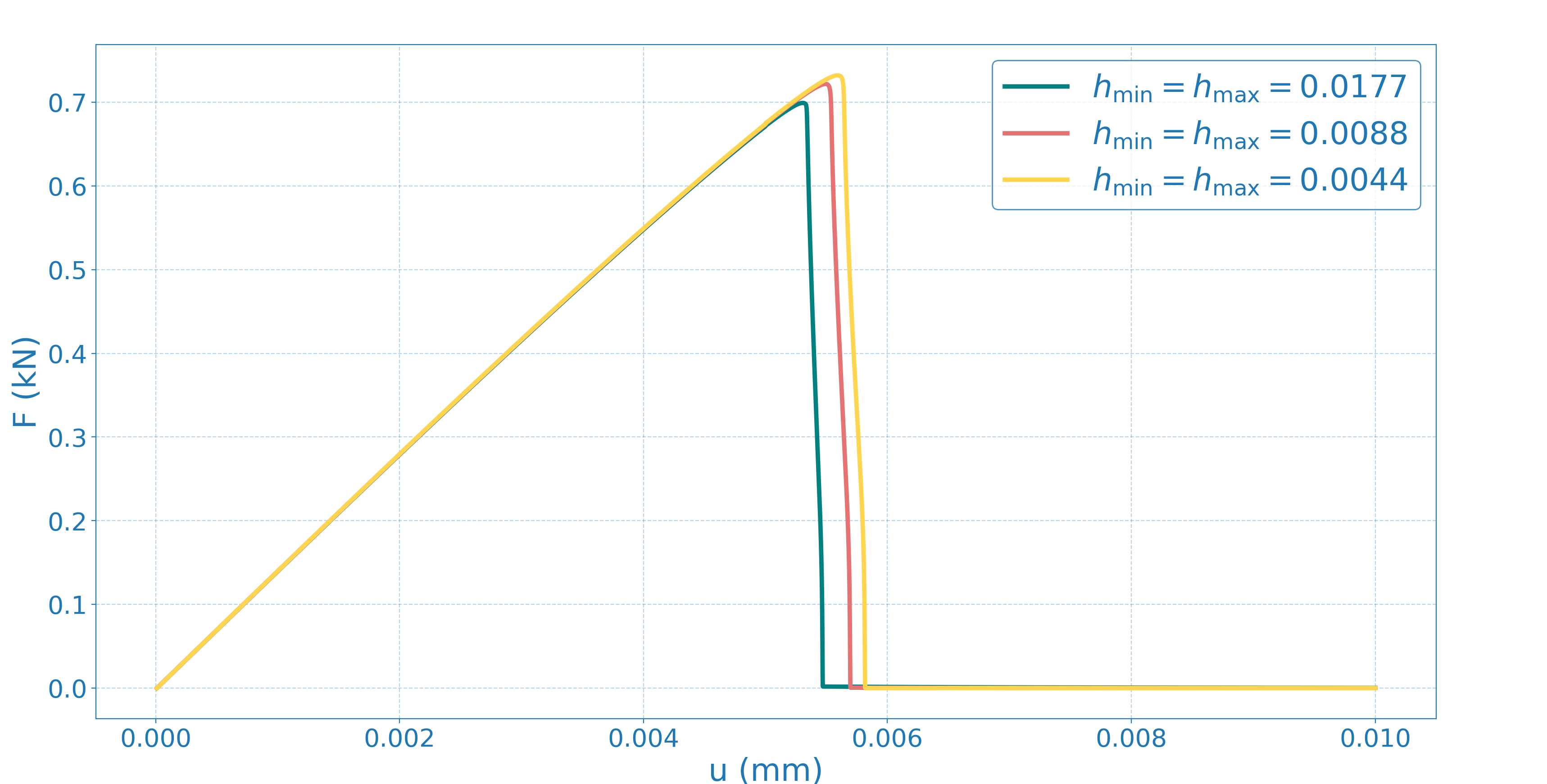}
  \caption{Load-displacement curves for the traction test with different mesh sizes using the hybrid–VD formulation.}
  \label{fig:modeI-load-displacement}
\end{figure}

A mesh convergence study was also performed.
Figure~\ref{fig:modeI-load-displacement} shows the load–displacement curves obtained with a family of uniform triangular meshes of size $h \in \{ \SI{0.0177}{\milli\metre},\SI{0.0088}{\milli\metre},\SI{0.0044}{\milli\metre}\}$ for the hybrid–VD formulation.
The horizontal axis represents the magnitude of the prescribed displacement $\vec{u}_{\rm left}$ at each time step applied on the left boundary, while the vertical axis reports the magnitude of the average reaction force computed on the same boundary.
A convergence trend can be observed as the mesh is refined, with the curves tending to stabilize for mesh sizes smaller than $h = \SI{0.01}{\milli\metre}$ .

\begin{figure}
  \centering
  \includegraphics[width=\textwidth]{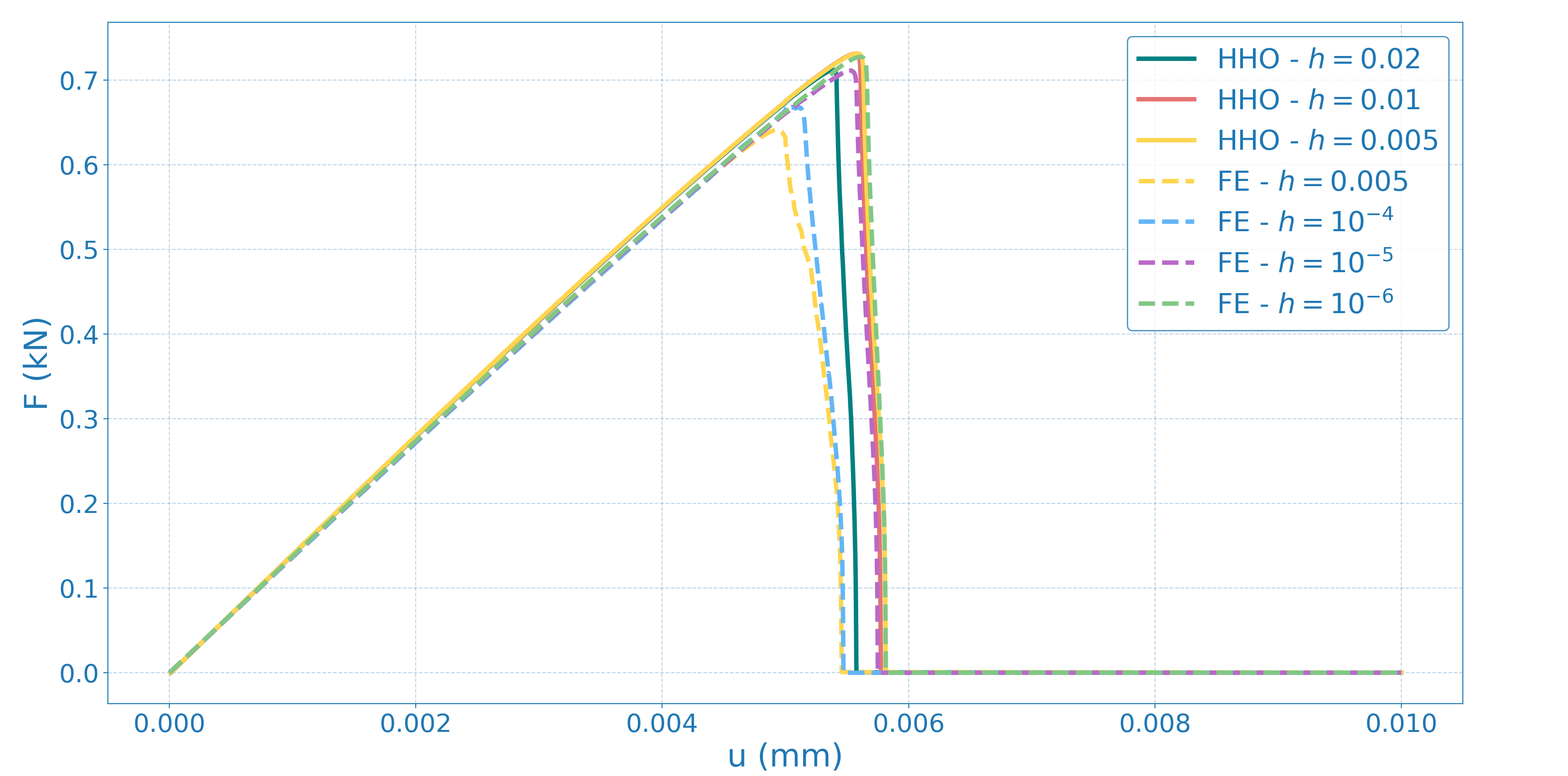}
  \caption{Load-displacement curves for the traction test with different mesh sizes using HHO (order $k=1$, $l=0$) and P2-P1 finite elements.}
  \label{fig:modeI-load-displacement-fem-hho}
\end{figure}

In Figure~\ref{fig:modeI-load-displacement-fem-hho} we compare the load-displacement curves obtained with our HHO scheme with those obtained using Lagrange P2-P1 finite elements on a family of triangular meshes.
The FEM simulations are performed with the FreeFEM++ software~\cite{Hecht:12}.
For HHO, we consider three uniform triangular meshes with sizes $h \in \{ \SI{0.02}{\milli\metre}, \SI{0.01}{\milli\metre}, \SI{0.005}{\milli\metre} \}$.
For FEM, we reuse the same uniform mesh with $h = \SI{0.005}{\milli\metre}$ used for HHO and, in addition, three non-uniform meshes generated adaptively with FreeFem++, featuring minimum element sizes respectively equal to $\SI{1e-4}{\milli\metre}$, $\SI{1e-5}{\milli\metre}$, and $\SI{1e-6}{\milli\metre}$.
We can observe that both methods tend to converge to the same solution as the mesh is refined, but with HHO showing a faster convergence despite the fact that the simulations are carried out on non-adaptive meshes.

\begin{figure}
  \centering
  \begin{subfigure}{0.45\textwidth}
    \includegraphics[height=6.5cm]{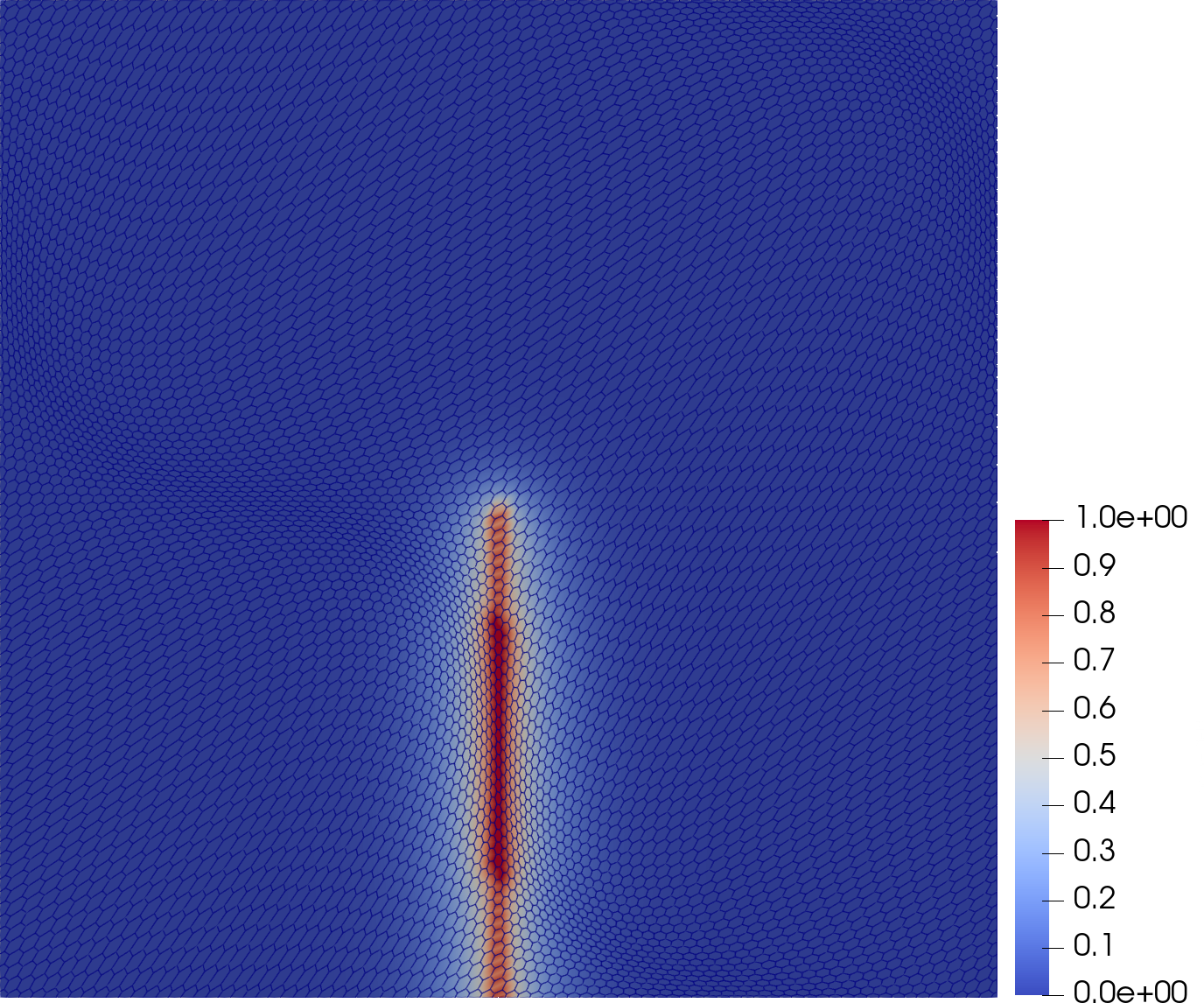}
    \caption*{$t=0$}
    \label{fig:hexa:0}
  \end{subfigure}
  \hfill
\begin{subfigure}{0.45\textwidth}
  \includegraphics[height=6.5cm]{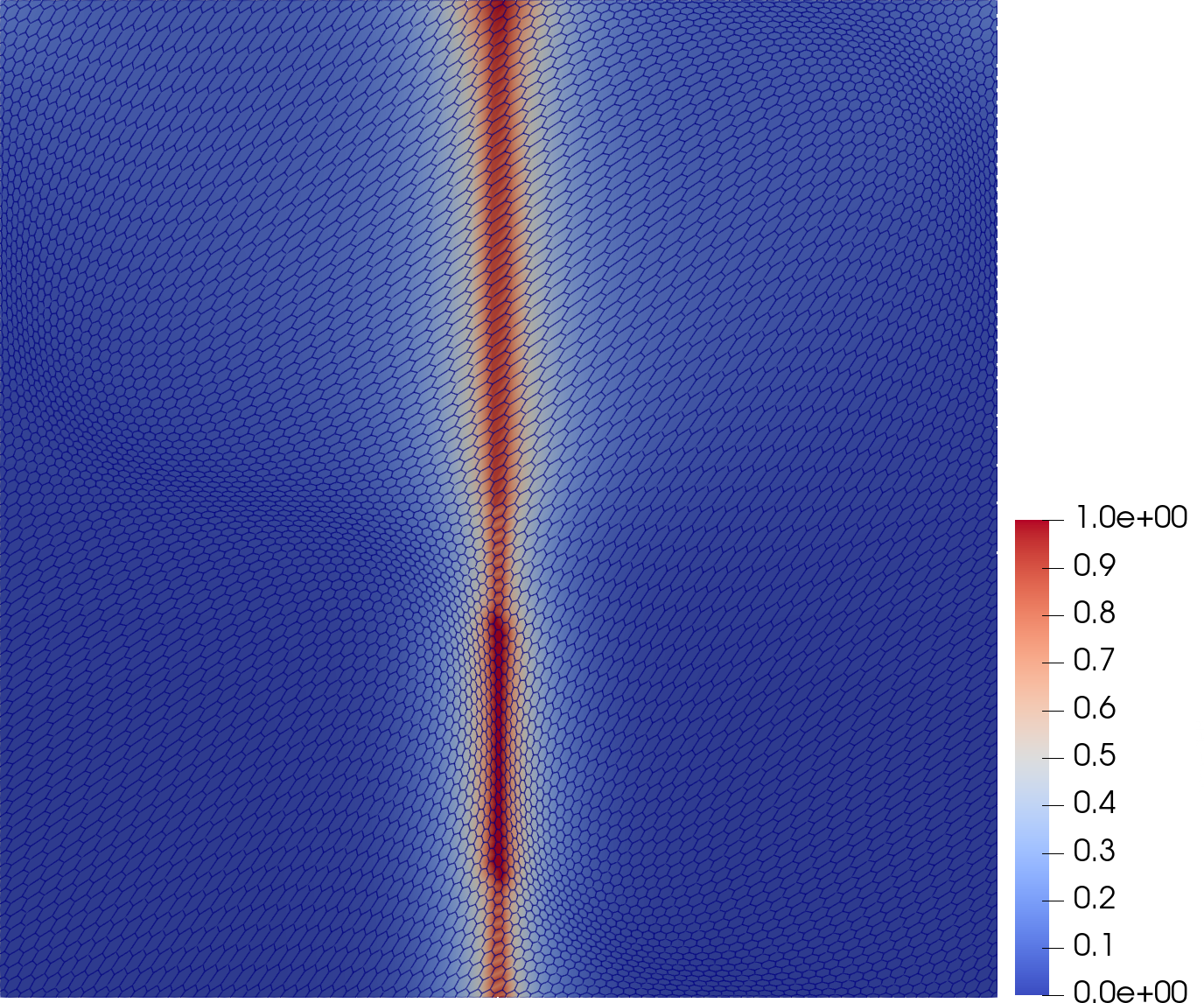}
  \caption*{$t=T_{end}$}
  \label{fig:hexa:end}
\end{subfigure}
\caption{Crack phase-field in mode I with HHO scheme on hexagonal meshes}
\label{fig:hexagonal-simulation}
\end{figure}

We also tested our scheme on a hexagonal mesh, which is shown in Figure~\ref{fig:hexagonal-simulation}.
The results are consistent with those obtained on triangular meshes, and confirm that the method behaves similarly across different mesh topologies.

\begin{figure}
  \centering
  \begin{subfigure}{0.45\textwidth}
    \includegraphics[height=6.5cm]{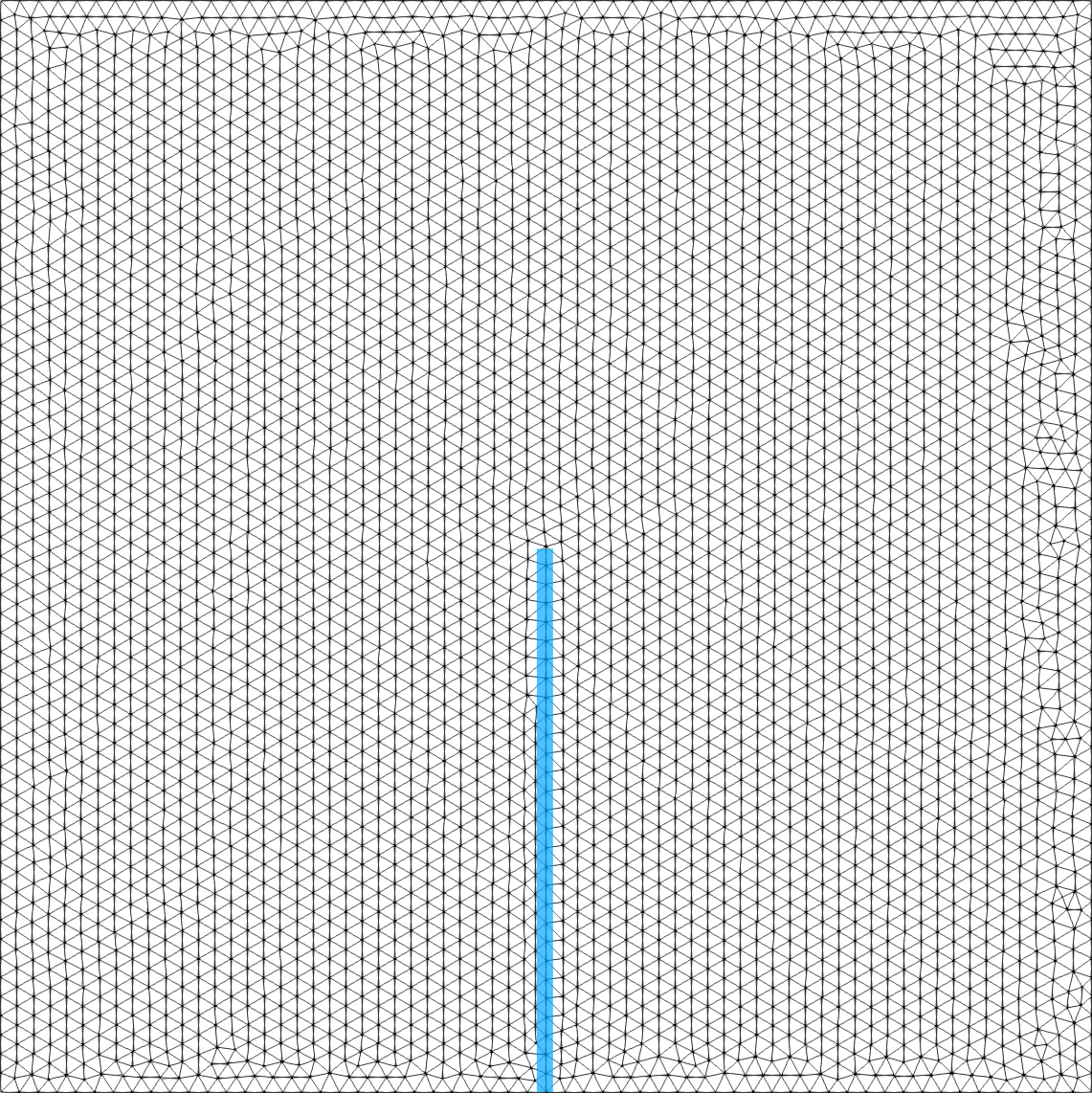}
    \caption{Initial notch represented by a physical discontinuity in the mesh, doubling the nodes along the notch line.}
    \label{fig:initial-notch-physical}
  \end{subfigure}
  \hfill
\begin{subfigure}{0.45\textwidth}
  \includegraphics[height=6.5cm]{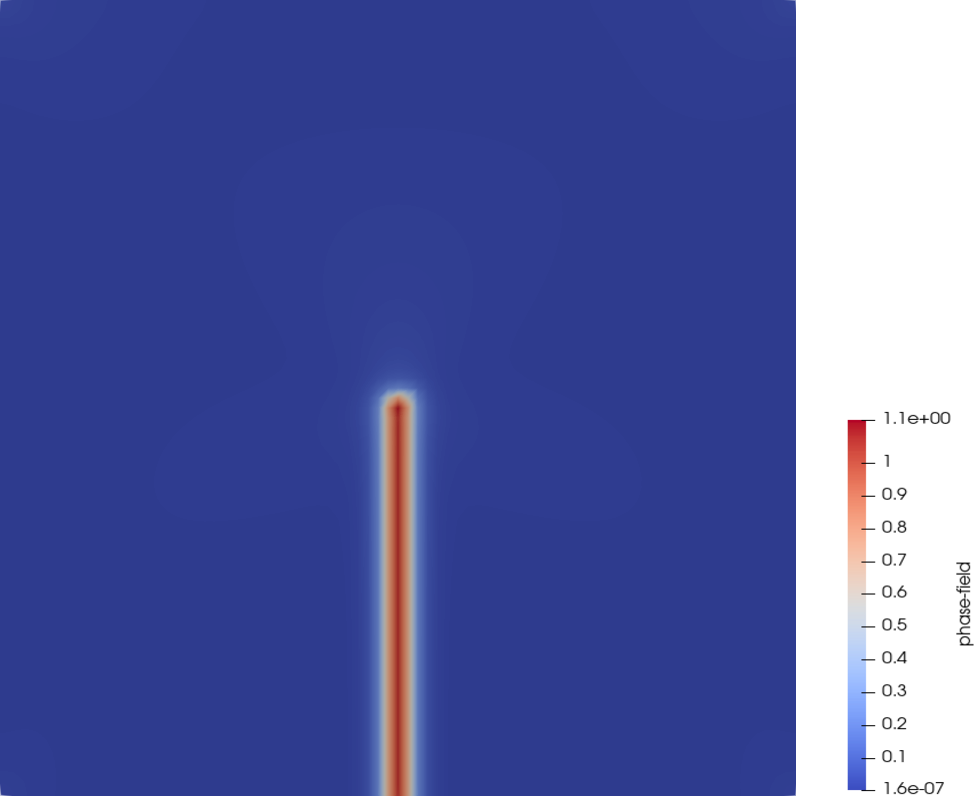}
  \caption{Initial notch represented by a high value of the initial history field $\calH_0$.\\}
  \label{fig:initial-notch-history}
\end{subfigure}
\caption{Two approaches to model the initial notch.}
\end{figure}

\begin{remark}[Initial notch]
  To model an initial notch in the domain, we considered two different approaches.
  In the tests with triangular meshes (see Figures~\ref{fig:phasefield-evolution-modeI}, \ref{fig:modeI-load-displacement}, \ref{fig:modeI-load-displacement-fem-hho}, \ref{fig:phasefield-evolution-modeII}, \ref{fig:modeII-load-displacement}, \ref{fig:modeII-load-displacement-fem-hho}),
  the discontinuity is modeled by introducing a physical cut in the mesh, obtained by duplicating the nodes along the notch line (highlighted in blue in Figure~\ref{fig:initial-notch-physical}).
  While this approach is straightforward, it requires the mesh to conform exactly to the notch geometry, which can be problematic for complex geometries.

  The second approach removes this problem by using an initial history field $\calH_0$ with high values along the notch line, thus effectively simulating a pre-existing crack without altering the mesh topology, as shown in Figure~\ref{fig:initial-notch-history}.
  This method offers greater flexibility, since the mesh no longer needs to match the notch geometry.
  It requires, however, careful calibration of the initial history field to correctly represent the notch. We used this approach for tests with polygonal meshes (see Figure~\ref{fig:hexagonal-simulation}).
\end{remark}

\subsection{Shear test}

\begin{figure}
  \centering
  \begin{subfigure}[t]{0.215\textwidth}
    \includegraphics[width=\textwidth]{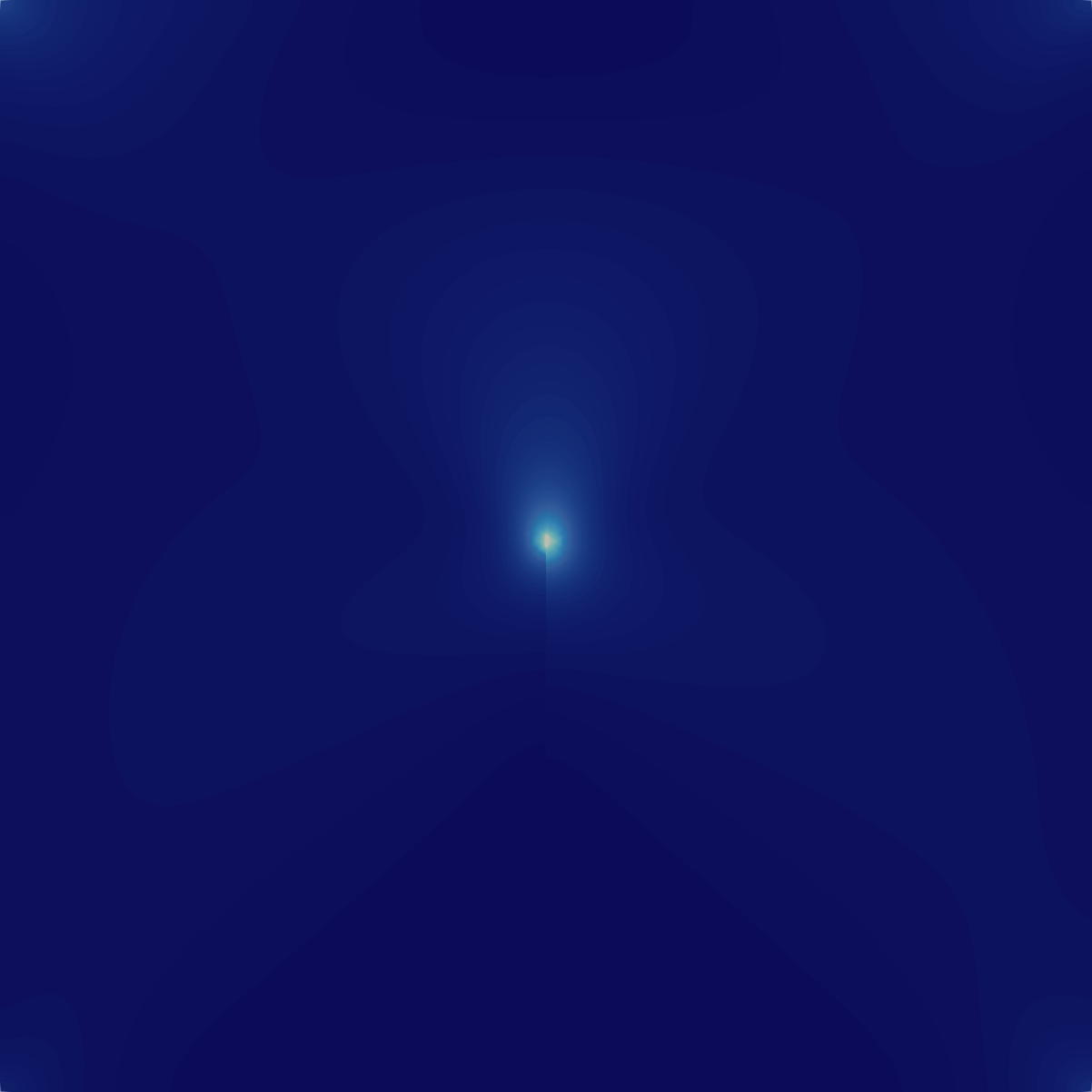}
  \end{subfigure}
  \hspace{0.0001cm}
  \begin{subfigure}[t]{0.215\textwidth}
    \includegraphics[width=\textwidth]{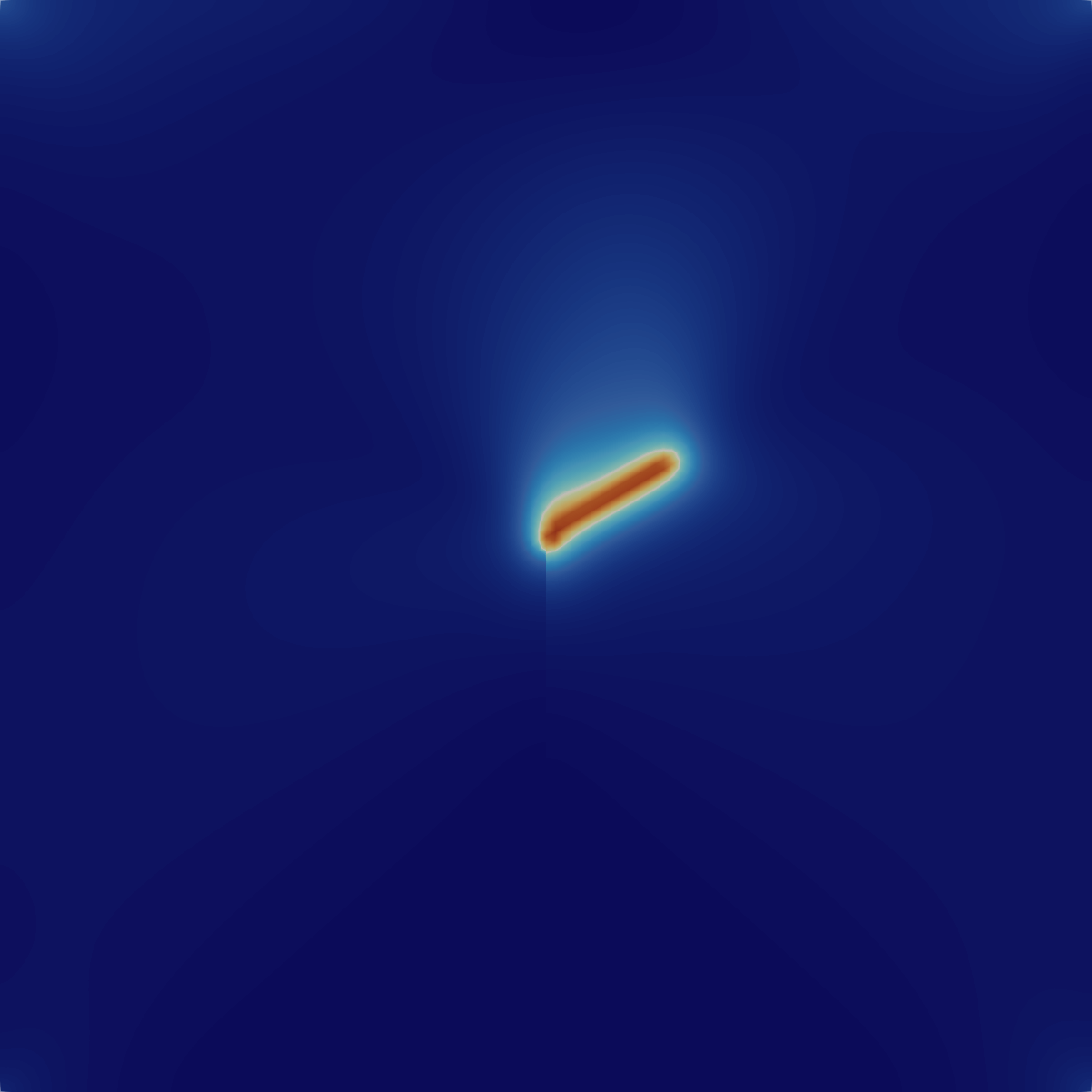}
  \end{subfigure}
  \hspace{0.0001cm}
  \begin{subfigure}[t]{0.215\textwidth}
    \includegraphics[width=\textwidth]{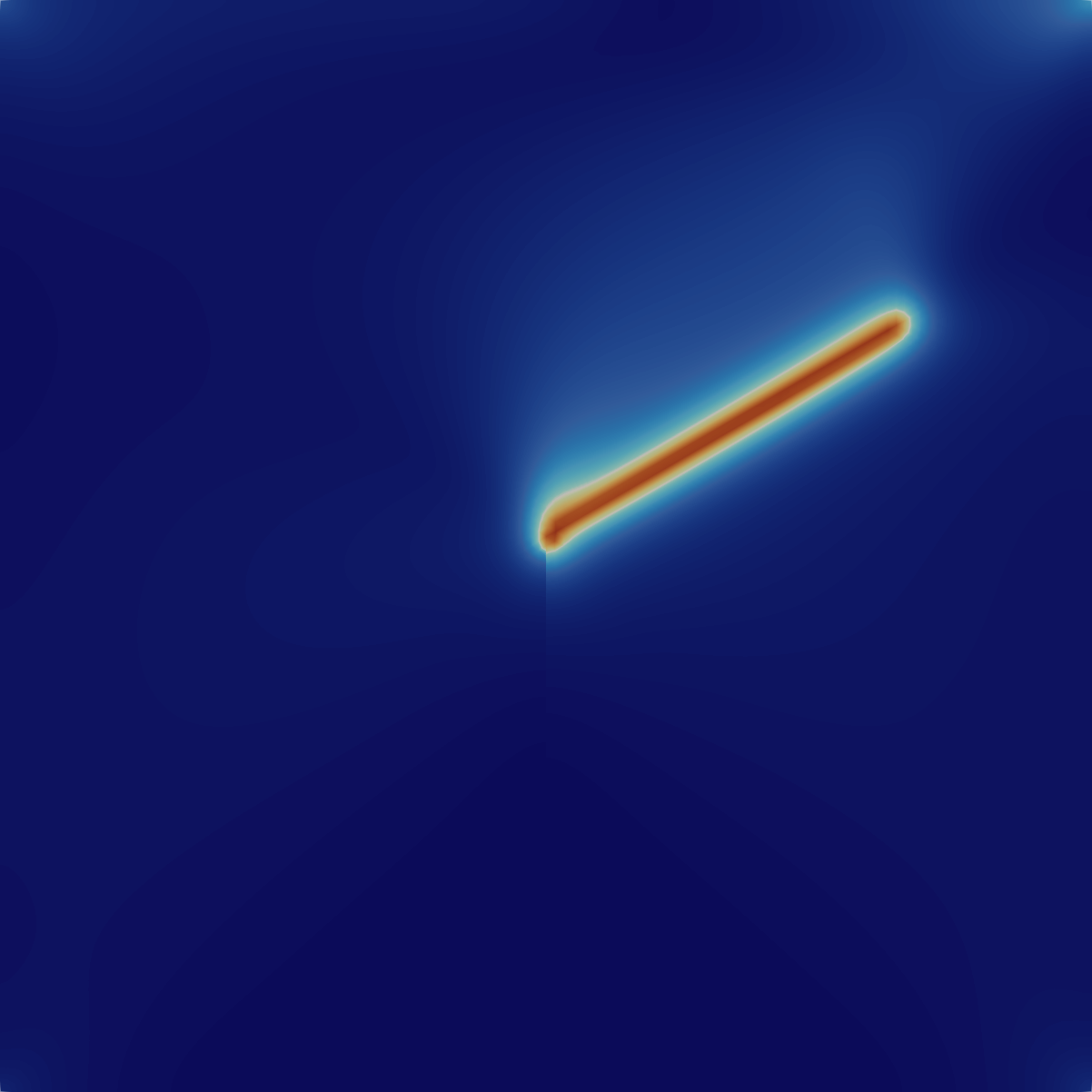}
  \end{subfigure}
  \hspace{0.0001cm}
  \begin{subfigure}[t]{0.215\textwidth}
    \includegraphics[width=\textwidth]{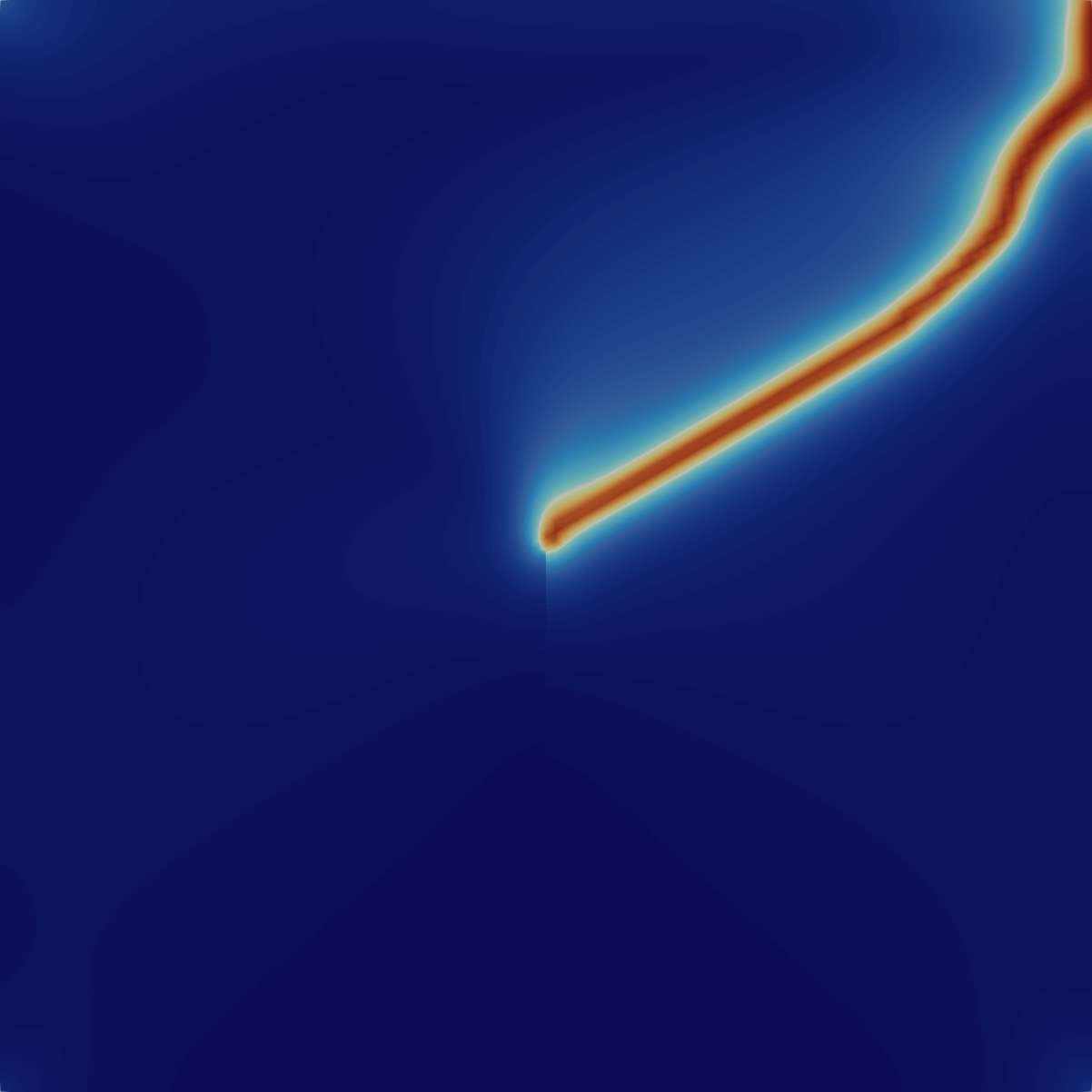}
  \end{subfigure}
  \hspace{0.001cm}
  \begin{subfigure}[t]{0.07\textwidth}
    \vspace{-3.2cm}
    \includegraphics[width=\textwidth]{Figures/hho_simulations/legend.png}
  \end{subfigure}

  \vspace{0.15cm}

  \begin{subfigure}[t]{0.215\textwidth}
    \includegraphics[width=\textwidth]{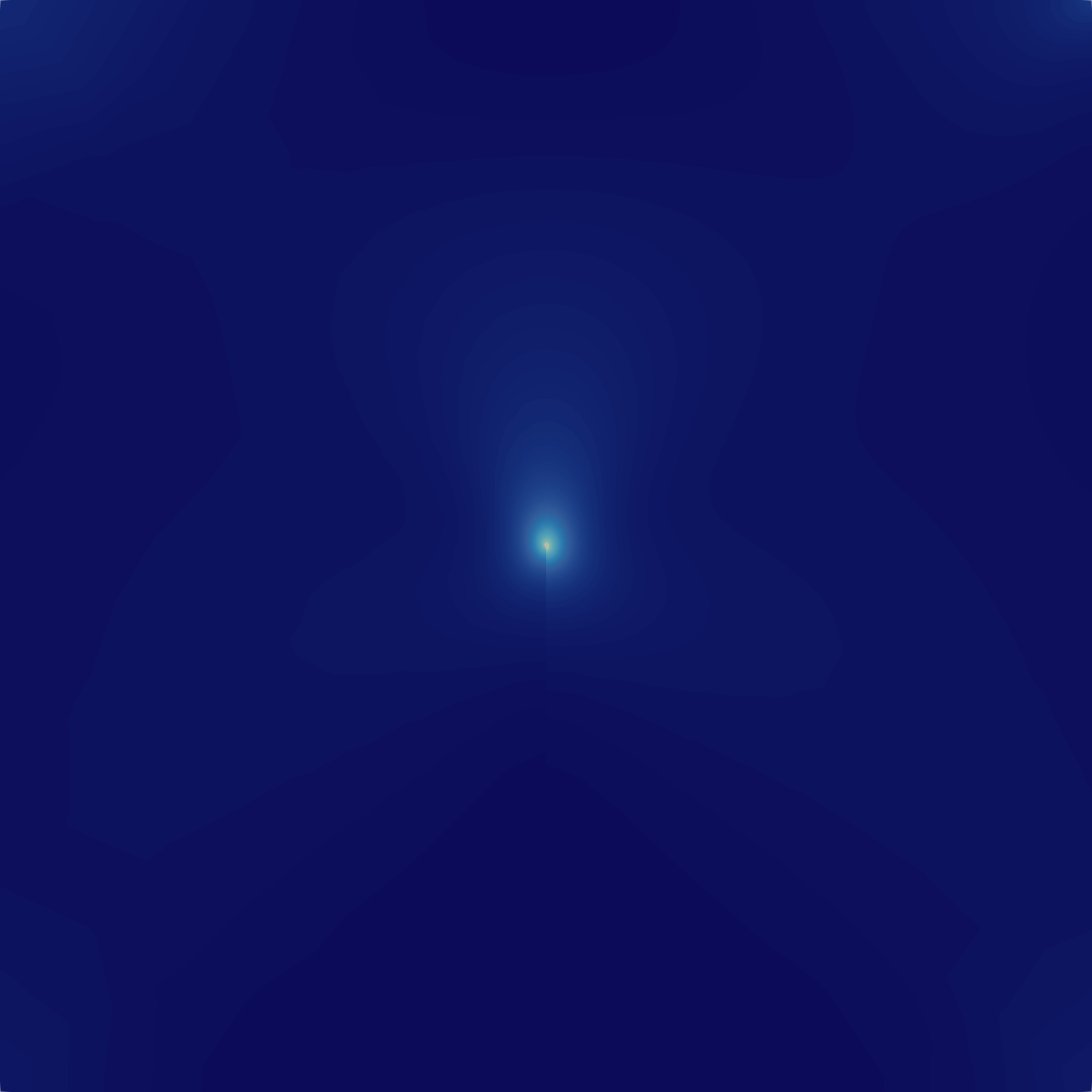}
    \centering{\footnotesize $t=70$}
  \end{subfigure}
  \hspace{0.0001cm}
  \begin{subfigure}[t]{0.215\textwidth}
    \includegraphics[width=\textwidth]{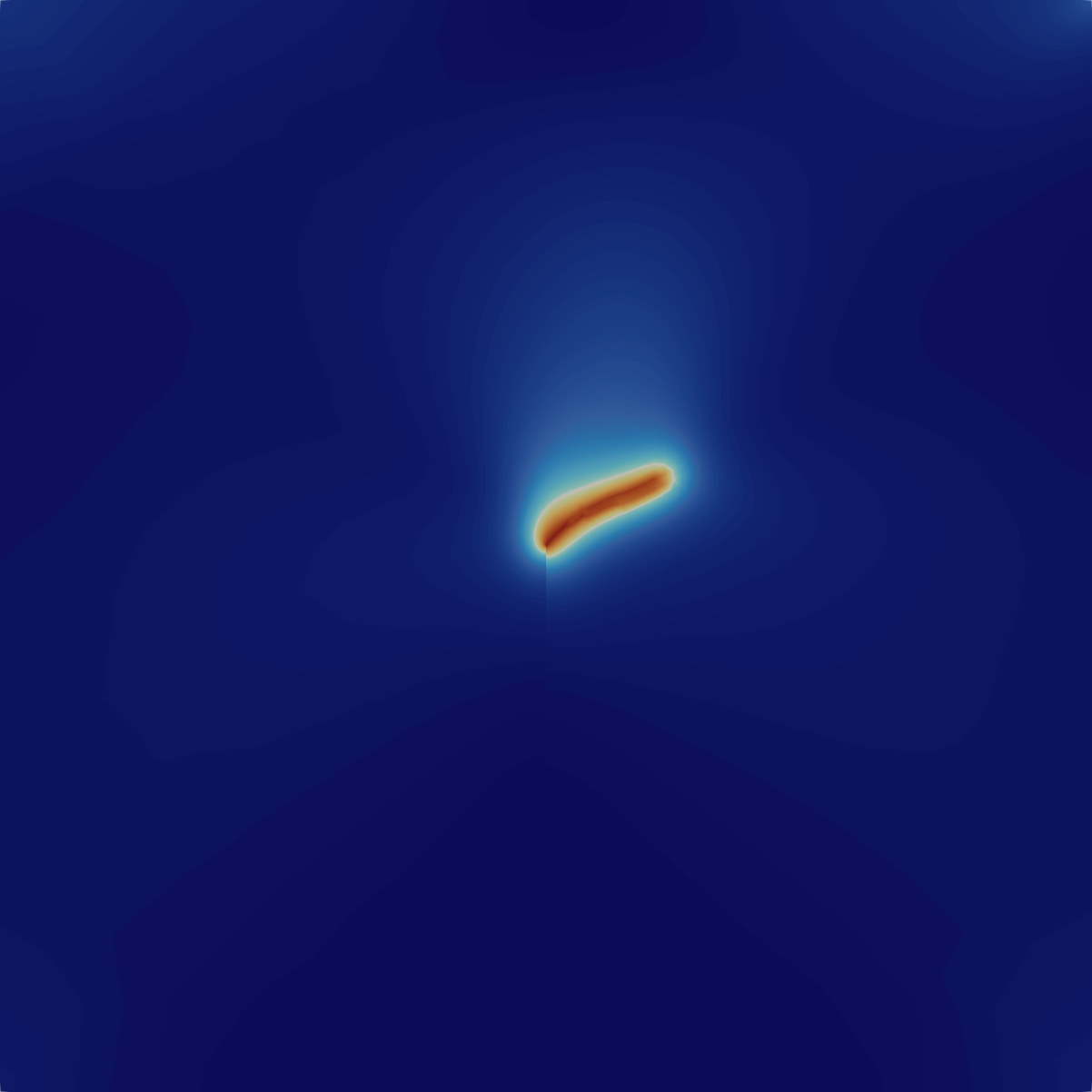}
    \centering{\footnotesize $t=100$}
  \end{subfigure}
  \hspace{0.0001cm}
  \begin{subfigure}[t]{0.215\textwidth}
    \includegraphics[width=\textwidth]{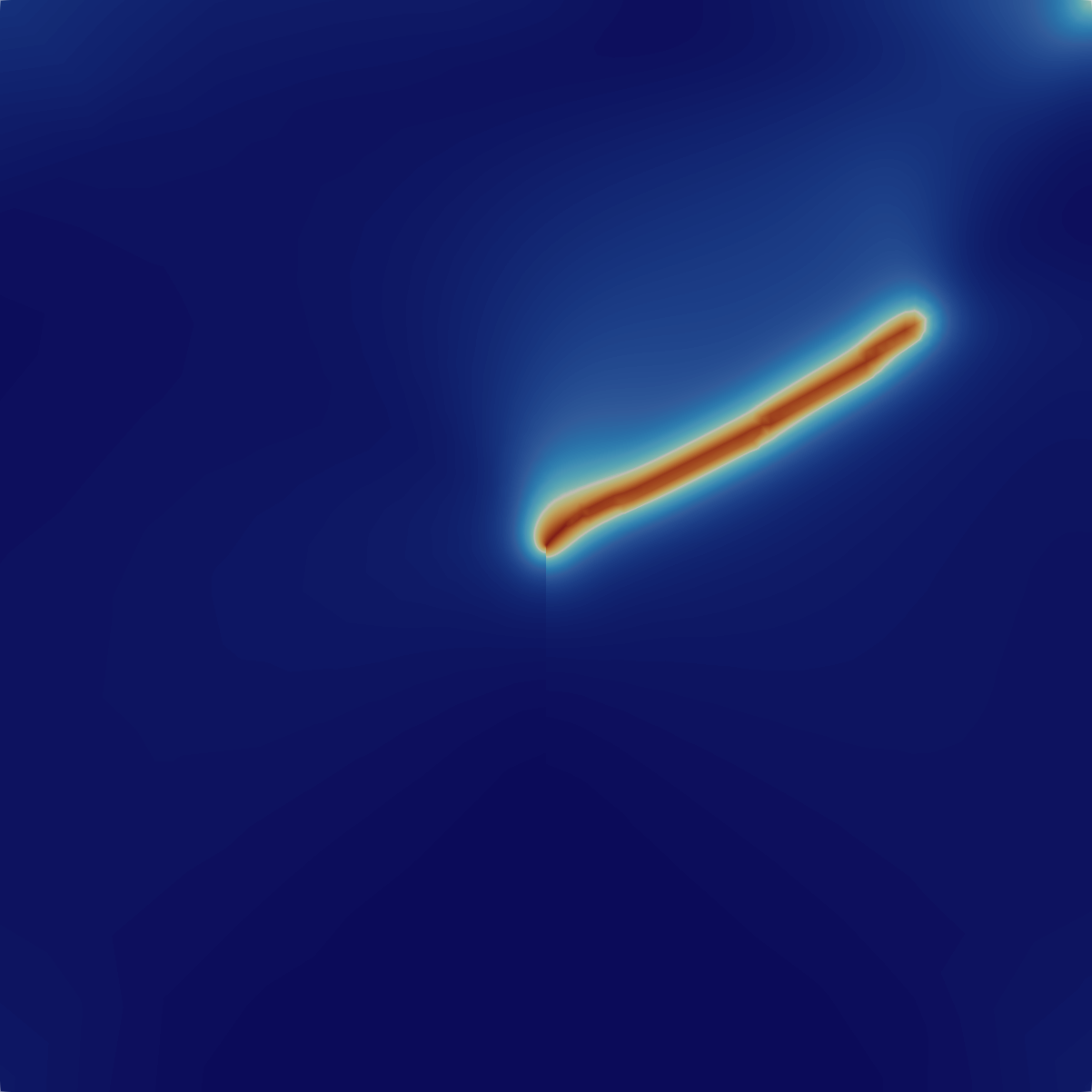}
    \centering{\footnotesize $t=130$}
  \end{subfigure}
  \hspace{0.0001cm}
  \begin{subfigure}[t]{0.22\textwidth}
    \includegraphics[width=\textwidth]{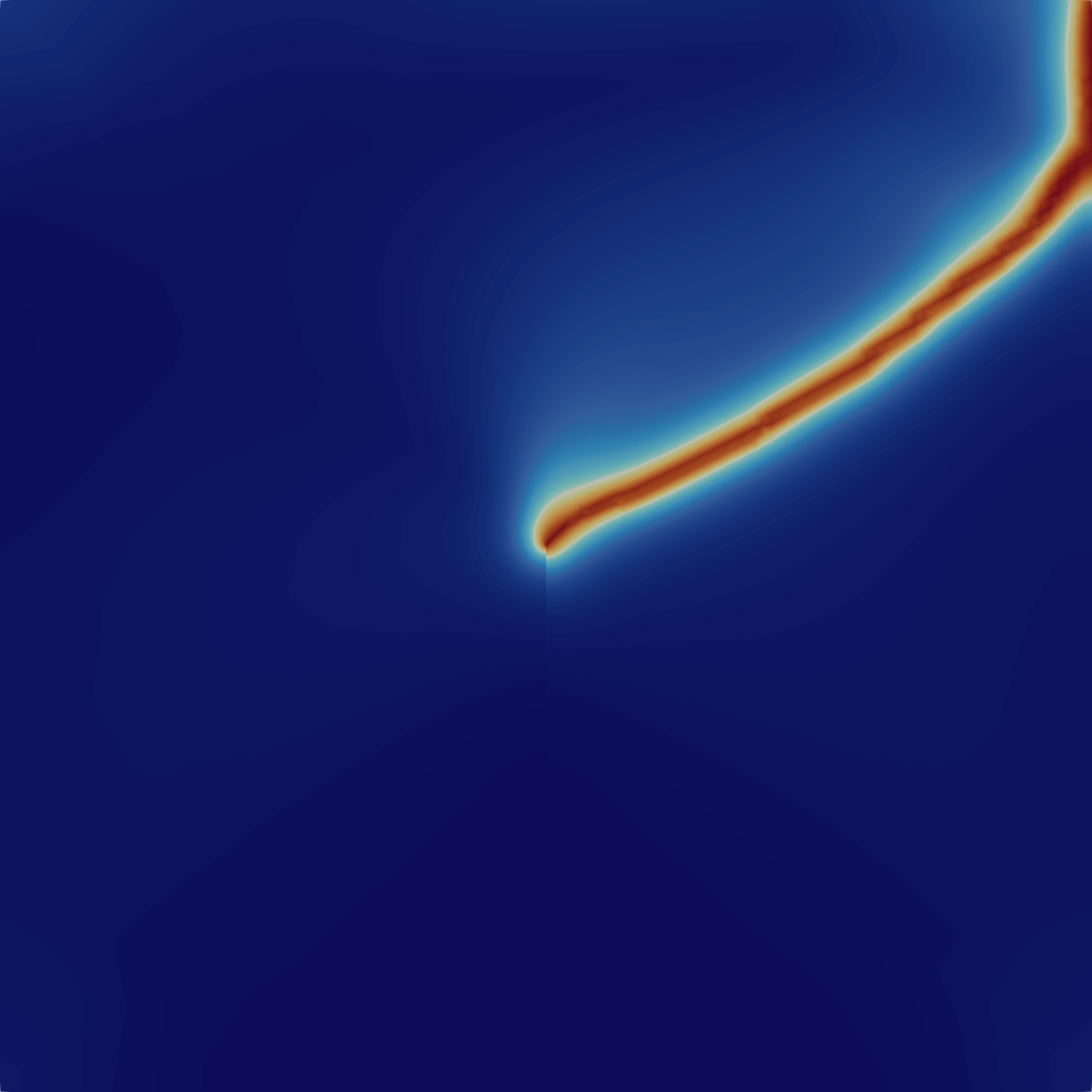}
    \caption*{\scriptsize $t=200$}
  \end{subfigure}
  \begin{subfigure}[t]{0.07\textwidth}
    \vspace{-3.2cm}
    \includegraphics[width=\textwidth]{Figures/hho_simulations/legend.png}
  \end{subfigure}

  \caption{Crack phase-field $\phi$ at different time steps for the mode II test.
  Top row: hybrid-VD formulation, bottom row: hybrid-SP energy decomposition.}
  \label{fig:phasefield-evolution-modeII}
\end{figure}

\begin{figure}
  \centering
  \includegraphics[width=0.4\textwidth]{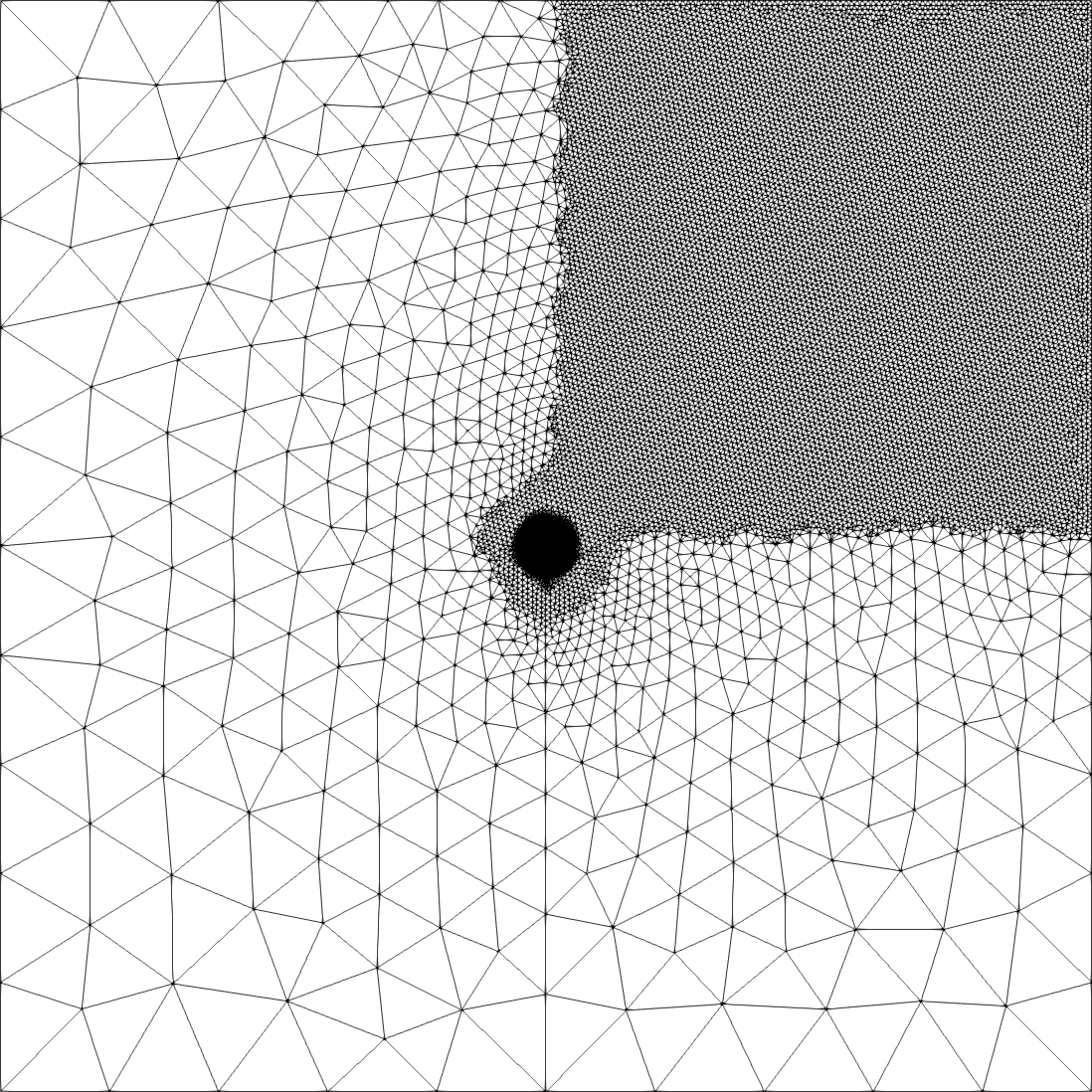}
  \caption{Locally refined mesh used for the shear test, with minimum element size $\SI{0.00813}{\milli\metre}$.}
  \label{fig:shear-mesh}
\end{figure}

The second benchmark corresponds to the \emph{mode II} shear loading test.
We consider the same domain as in the traction test, but with a displacement applied on the left boundary in the upward direction, as shown on the right in Figure~\ref{fig:domain}.
The prescribed displacement increases monotonically, starting from zero and incremented linearly by $\SI{1e-5}{\milli\metre}$ at each time step.
We use triangular meshes locally refined where the crack is expected to propagate (see Figure~\ref{fig:shear-mesh}).
The Lamé parameters $\lambda$ and $\mu$, the critical energy release rate $G_c$, and the regularization length $\ell$ are the same as in the traction test.
Also in this case, several simulations are performed to investigate the convergence of the method and to compare the results obtained with the HHO scheme to those reported in the literature using other discretization methods.
In Figure~\ref{fig:phasefield-evolution-modeII} we show the evolution of $\phi$ at different time steps for hybrid-VD and hybrid-SP formulations.

\begin{figure}
  \centering
  \includegraphics[width=\textwidth]{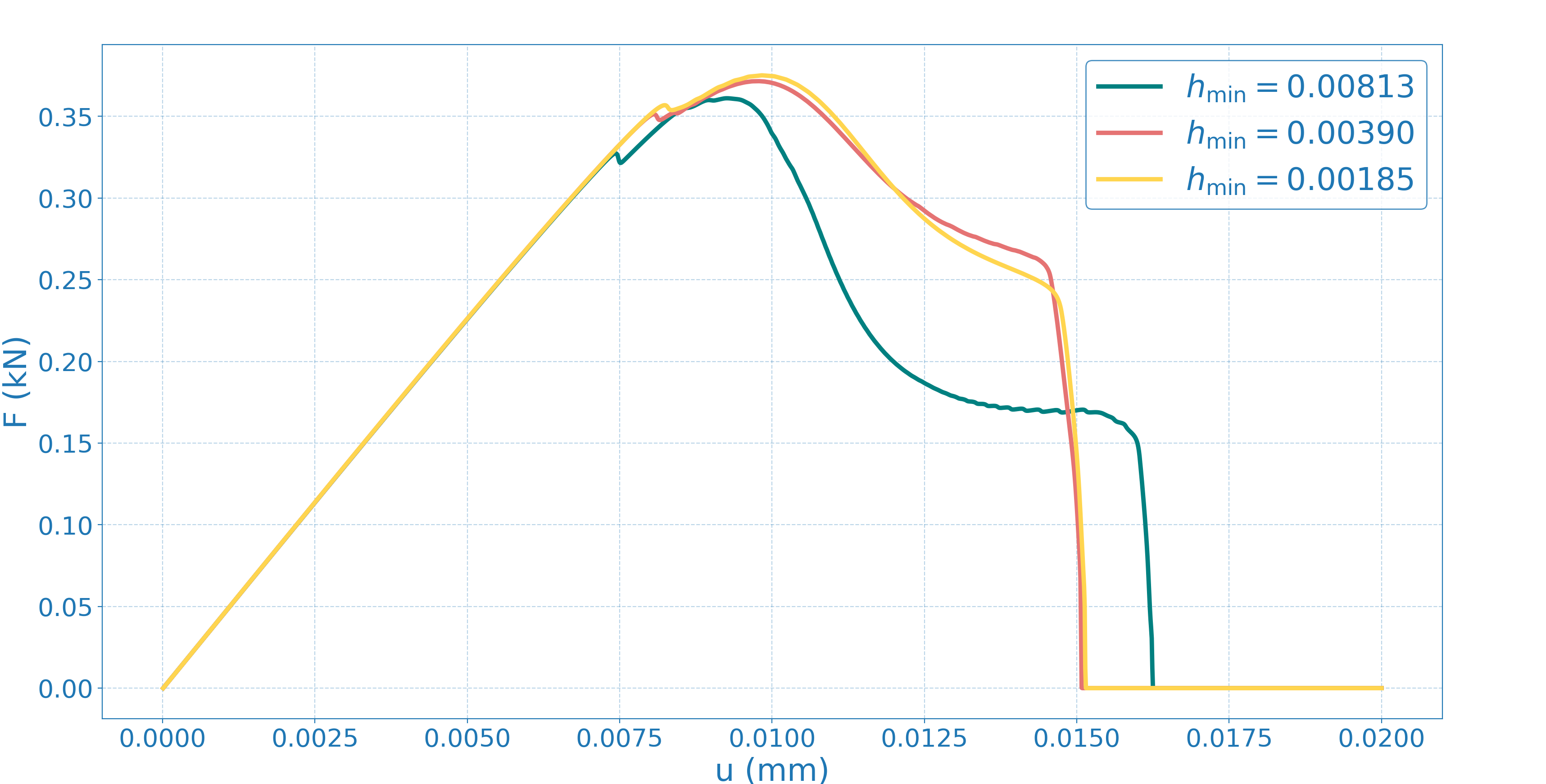}
  \caption{Load-displacement curves for the shear test with different mesh sizes using the hybrid–VD formulation.}
  \label{fig:modeII-load-displacement}
\end{figure}

The load–displacement curves for the hybrid–VD formulation obtained with a family of locally refined meshes similar to the one depicted in Figure~\ref{fig:shear-mesh} with minimum element sizes in $\{ \SI{0.00813}{\milli\metre}, \SI{0.00390}{\milli\metre}, \SI{0.00185}{\milli\metre} \}$, are shown in Figure~\ref{fig:modeII-load-displacement}.
As in the traction test, the horizontal axis represents the magnitude of the prescribed displacement $\boldsymbol{u}_{\rm left}$ at each loading increment, applied on the left boundary, while the vertical axis reports the magnitude of the average reaction force computed on the same boundary.
A convergence trend can be observed as the mesh is refined, with the curves tending to stabilize for mesh sizes smaller than $\SI{0.004}{\milli\metre}$.
Notice that, for the coarsest mesh, the mesh size is larger than $\ell$, which can justify the significant difference with respect to the two finest meshes.

\begin{figure}
  \centering
  \includegraphics[width=\textwidth]{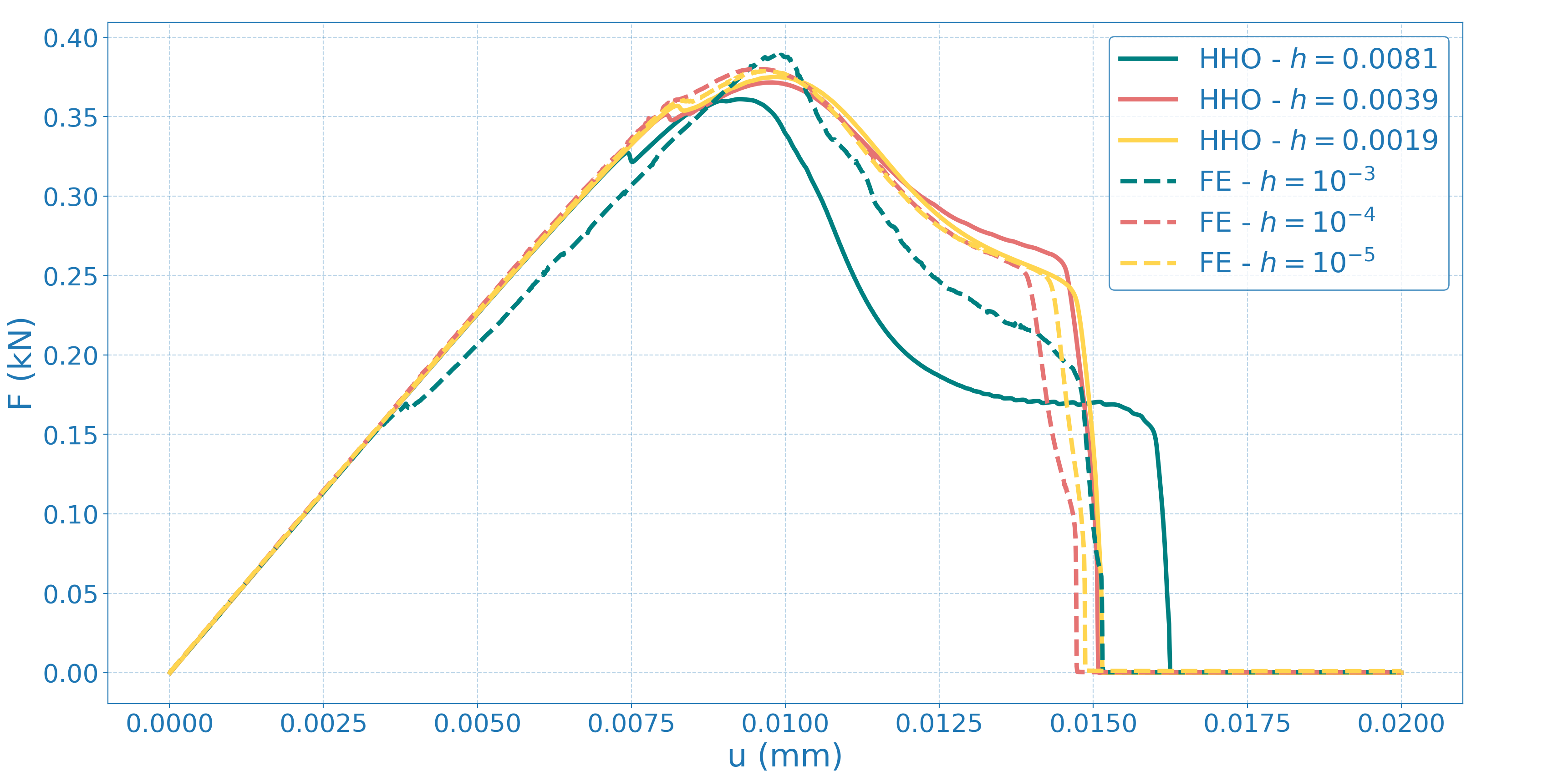}
  \caption{Load-displacement curves for the shear test with different mesh sizes using HHO (order $k=1$, $l=0$) and P2-P1 finite elements.}
  \label{fig:modeII-load-displacement-fem-hho}
\end{figure}

In Figure~\ref{fig:modeII-load-displacement-fem-hho} we provide a comparison with Lagrange P2-P1 FEM discretization on a family of locally refined triangular meshes. For HHO we consider the same meshes as for the previous convergence test, similar to the one shown in Figure~\ref{fig:shear-mesh}, while for FEM we use three non-uniform meshes generated adaptively with FreeFem++, featuring minimum element sizes in $\{ \SI{1e-3}{\milli\metre}, \SI{1e-4}{\milli\metre}, \SI{1e-5}{\milli\metre} \}$.
We can observe that both methods tend to converge to the same solution as the mesh is refined and, once again, HHO appears to reach the exact solution faster than FEM.


\bibliographystyle{abbrv}
\bibliography{hho-phasefield}

\end{document}